\documentclass{article}

\usepackage{amsmath}
\usepackage{camnum}
\usepackage{graphicx}
\usepackage{harvard}
\usepackage{xcolor}
\usepackage[top=2.5cm, bottom=2.5cm, left=2.5cm, right=2.5cm]{geometry}

\newtheorem{assumption}{Assumption}
\newtheorem{definition}{Definition}
\newtheorem{example}{Example}
\newtheorem{remark}{Remark}

\newcommand{\ee}{{\mathrm e}}
\newcommand{\ii}{{\mathrm i}}

\definecolor{darkgreen}{RGB}{0,100,0}

\begin{document}

\title{Symmetric Zassenhaus splittings\\
	for bounded and unbounded operators} 

\author{
Paul Bergold
\footnote{Institute of Mathematics,
	University of Gda\'nsk, Gda\'nsk, Poland.
	Corresponding author: \texttt{paul.bergold@ug.edu.pl}}
\and
Arieh Iserles
\footnote{
	Department of Applied Mathematics and Theoretical Physics,
	University of Cambridge, Cambridge, United Kingdom}
\and
Karolina Kropielnicka
\footnote{	Institute of Mathematics,
	University of Gda\'nsk, Gda\'nsk, Poland}
}

\thispagestyle{empty}
\maketitle
\begin{abstract}
	We establish rigorous local-error bounds for the entire family of symmetric Zassenhaus splittings for possibly unbounded operators, under explicit assumptions on commutator bounds, regularity-preserving flows, and conjugation identities.
	Truncated symmetric Baker--Campbell--Hausdorff expansions yield an explicit coefficient formula for the Zassenhaus exponents and identify vanishing commutator expressions.
	In order to extend these cancellations to the realm of unbounded operators, we prove a general theorem on finite expansions of multiple operator conjugations with controlled integral remainders.
	For splittings of order $2n$, our argument gives the sufficient commutator-length bound $4n-1$, revealing regularity requirements that are not apparent from the formal order calculation.
	Using the Schr\"odinger equation with a quadratic potential as an example, we show that the order-dependent regularity requirements can reduce to one fixed weighted Sobolev space for all splitting orders.
	
\end{abstract}

\tableofcontents
\section{Introduction}
Splittings are a powerful methodology in the numerical solution of partial differential equations.
Their main idea is to approximate the evolution of a complicated problem by composing the flows of simpler subproblems.
Consider, for instance, the linear evolution equation
\begin{equation*}
	u'(t)
	=(A+B)u(t),\qquad
	u(0)
	=u_0,
\end{equation*}
where $A$ and $B$ are possibly unbounded linear operators on a Hilbert space $\mathcal{H}$.
Assume that $A$, $B$, and a suitable realization of $A+B$ generate strongly continuous semigroups, denoted by $\ee^{tA}$, $\ee^{tB}$, and $U(t)$, respectively.
Direct approximation of $U(t)$ may be expensive, even when the subproblems $u'=Au$ and $u'=Bu$ can be solved efficiently.
Splitting methods exploit this structure by advancing the numerical solution from $n\tau$ to $(n+1)\tau$, with time step $\tau>0$, through a composition of these simpler flows.

For example, under suitable assumptions the \emph{Lie--Trotter splitting} $u^{n+1}=\ee^{\tau A}\ee^{\tau B}u^n$ and the \emph{Strang splitting} $u^{n+1}=\ee^{\frac12\tau A}\ee^{\tau B}\ee^{\frac12\tau A}u^n$ have orders one and two, respectively.
There exists an extensive theory of splitting methods, their advantages and limitations; see \cite{Blanes2024,McLachlan2002} and references therein.
Higher-order schemes may be constructed and analysed by several approaches, including matching finite Taylor expansions and exploiting the Baker--Campbell--Hausdorff (BCH) formula.

An important example, which also motivates the present work, is the linear Schr\"odinger equation with a time-independent potential,
\begin{equation*}
	\ii\partial_t u
	=-\Delta u+V(x)u,
\end{equation*}
corresponding to the choice $A=\ii\Delta$ and $B=-\ii V$ on $\mathcal{H}=L^2(\BB R^d)$.
The kinetic and potential subflows can often be computed efficiently, even when their generators are unbounded, as is the case for the Laplacian and multiplication by an unbounded potential, such as the Coulomb potential.
The domain and regularity issues arising from these unbounded operators are therefore relevant to practical applications and are not merely theoretical.
The formal order of a splitting need not be attained if the required domain and regularity assumptions are not satisfied: in particular settings, Strang splitting can exhibit order reduction, with only first-order convergence instead of the anticipated second order \cite{Jahnke2000,Iserles2024,Kropielnicka2025}.
Coulomb interactions arise in computational models of atomic dynamics \cite{Marante2014}; related work on hydrogenic Hamiltonians illustrates that product-formula convergence may be substantially slower for certain states \cite{Burgarth2024}.
Establishing the accuracy of splitting methods for unbounded operators is therefore important both for mathematical analysis and practical computation.

For bounded $A$ and $B$, operator power series converge in norm, and finite Taylor expansions or suitably truncated BCH constructions can be justified without domain restrictions.
This argument cannot simply be carried over to unbounded generators.
Infinite operator expansions generally need not converge in operator norm, and even individual products and commutators may fail to act on a common domain.
Finite Taylor expansions \emph{with integral remainders} can nevertheless be rigorously justified for sufficiently regular data.
The difficulty is then to control the derivatives and remainders: direct expansion of the exact evolution and of all splitting factors can entail high powers and mixed products of unbounded operators, while expansions of their conjugations (\emph{i.e.}, terms of the form $\ee^{tX}Y\ee^{-tX}$) generate nested commutators whose well-definedness must be established.
The question is not whether formal cancellations can be written down, but under what operator and regularity assumptions they imply an error estimate in the underlying Hilbert space $\mathcal{H}$.\\

The subject of this paper is the \emph{symmetric Zassenhaus splitting} (sZS), originally introduced as an efficient approach to the linear Schr\"odinger equation in the semiclassical regime \cite{Bader2014}.
It constructs a higher-order approximation by inserting exponentials of Lie polynomials into a symmetric composition of the flows generated by $A$ and $B$.
The entire family can be written as
\begin{equation*}
	Z_{2n-1}(\tau)
	=\ee^{\frac12\tau A}\ee^{\frac12\tau B}\left(\prod_{\ell=2}^{n}\ee^{\mathcal{R}_{2\ell-1}(\tau)}\right) \left(\prod_{\ell=2}^{n}\ee^{\mathcal{R}_{2(n-\ell)+3}(\tau)}\right)\ee^{\frac12\tau B}\ee^{\frac12\tau A},
\end{equation*}
where $\mathcal{R}_{2\ell-1}(\tau)=\tau^{2\ell-1}\mathcal{R}_{2\ell-1}(1)$ and each $\mathcal{R}_{2\ell-1}(1)$ is a linear combination of nested commutators containing exactly $2\ell-1$ occurrences of $A$ and $B$.
The case $n=1$ recovers the Strang splitting, while successive commutator corrections produce schemes of arbitrarily high (necessarily even) order.
In particular, for Schr\"odinger operators these commutators often simplify, making their exponentials especially attractive in semiclassical calculations \cite{Bader2014}.
Arnal, Casas and Chiralt developed a recursive construction of the symmetric Zassenhaus exponents and proved convergence of the infinite factorization in a Banach-algebra setting \cite{Arnal2017}.
Such results do not, by themselves, justify the asserted orders for unbounded operators.

There is an extensive theory of splitting methods for unbounded operators, particularly for compositions of the flows generated by $A$ and $B$ \cite{Hansen2009}; recent work also provides commutator-based error estimates for unitary product formulas \cite{Casas2026}.
Zassenhaus splittings introduce an additional layer: the correction terms themselves must generate suitable flows, and their repeated conjugations must be controlled on compatible domains.
An important contribution is the fourth-order analysis by Lasser and Lubich for the semiclassical Schr\"odinger equation \cite{Lasser2020}.
Their proof combines a symmetrized defect, integral representations of conjugated operators, and algebraic cancellations to derive a global fourth-order error estimate.
It is, however, tailored to the fourth-order Zassenhaus formula and does not directly provide a systematic analysis of the hierarchy of commutator-generated flows at higher orders.
It is this missing arbitrary-order framework that we develop here.

Our main result establishes rigorous local-error estimates for symmetric Zassenhaus splittings of arbitrarily high order for possibly unbounded operators.
More precisely, for every fixed integer $n\ge1$, subject to the stated conditions on commutator bounds, regularity-preserving flows, and integral conjugation identities, we show that
\begin{equation}\label{eq:intro-local-error-bound}
	\|U(\tau)\psi-Z_{2n-1}(\tau)\psi\|_{\mathcal{H}}
	\le c\tau^{2n+1}\|\psi\|_{\mathcal{D}},\quad
	\psi\in\mathcal{D},\quad
	0\le\tau\le\tau_0.	
\end{equation}
Here $\mathcal{D}$ is a Banach space with norm $\|\cdot\|_{\mathcal{D}}$, densely and continuously embedded in $\mathcal{H}$, preserved by the relevant flows, and chosen so that the required commutators are bounded from $\mathcal{D}$ to $\mathcal{H}$.
The precise assumptions are stated in Theorem~\ref{thm:unbounded-local-error}.

The proof combines a bounded-operator analysis with finite expansions that remain valid for unbounded operators.
For bounded operators, truncated symmetric BCH expansions yield an explicit coefficient formula for the Zassenhaus exponents (Theorem~\ref{thm:zassenhaus-exponents}) and identify the leading local-error contribution (Theorem~\ref{thm:bounded-local-error}).
To transfer the resulting error order to unbounded operators, we follow the symmetrized-defect approach of Lasser and Lubich \cite{Lasser2020}.
This expresses the local error in terms of a defect whose lower-order coefficients are combinations of nested commutators.
The bounded-operator analysis shows that these combinations vanish identically.
Since this vanishing is algebraic, it also holds for unbounded operators wherever the required products are defined.
More precisely, for $E_{2n-1}(\tau)=U(\tau)-Z_{2n-1}(\tau)$, the symmetrized-defect approach gives, on a suitable common domain,
\begin{equation}\label{eq:intro-duhamel}
	E_{2n-1}(\tau)
	=\int_0^\tau U\!\left(\frac{\tau-s}{2}\right)N_{2n-1}(s)\,U\!\left(\frac{\tau-s}{2}\right)\,\mathrm{d}s.
\end{equation}
The key step is to factor the defect as $N_{2n-1}(\tau)=G_{2n-1}(\tau)Z_{2n-1}(\tau)$.
For bounded operators, the Taylor coefficients $K_k$ of $G_{2n-1}$ are homogeneous Lie polynomials, and we show that $K_0=\cdots=K_{2n-1}=0$ (Theorem~\ref{thm:vanishing-defect-coefficients}).
These coefficients vanish algebraically, independently of the particular operators wherever the required products are defined.
The main analytic task is therefore to justify the corresponding finite expansion of $G_{2n-1}$ for unbounded operators with a controlled remainder.
The vanishing coefficients then give $G_{2n-1}(\tau)=\mathcal{O}(\tau^{2n})$.
Together with local boundedness of $Z_{2n-1}$ on $\mathcal{D}$, the factorization and \eqref{eq:intro-duhamel} yield the local-error bound \eqref{eq:intro-local-error-bound}.

In the course of proving the main theorem, we develop a general analytic tool of independent interest: Theorem~\ref{thm:multiple-conjugation-expansion} establishes finite expansions with integral remainder bounds in the space $\mathcal{L}(\mathcal{D},\mathcal{H})$ of bounded linear operators from $\mathcal{D}$ to $\mathcal{H}$, for multiple conjugations of the form
\begin{equation*}
	\left(\ee^{t^{k_m}X_m}\cdots\ee^{t^{k_1}X_1}\right)Y\left(\ee^{-t^{k_1}X_1}\cdots\ee^{-t^{k_m}X_m}\right),
\end{equation*}
where $X_1,\dots,X_m$ and $Y$ are homogeneous Lie polynomials in $A$ and $B$.
The result identifies not only the coefficients of the finite expansion but also a concrete, sufficient class of commutators controlling its remainder.
Unlike a formal sBCH calculation, this argument makes the domain and regularity requirements explicit.
In particular, the conjugation theorem is not specific to the Zassenhaus arrangement of factors.
Together with the algebraic cancellation mechanism, it suggests a transferable proof strategy for other exponential splittings whose factors involve commutator-generated flows, although such extensions are not proved here.

The remainder bounds in Theorem~\ref{thm:multiple-conjugation-expansion} reveal a further aspect of the analysis: regularity may be governed by commutators that do not appear in the formally truncated Zassenhaus expansion.
For the symmetric Zassenhaus splitting $Z_{2n-1}$, our proof requires control of commutators of length at most $4n-1$, whereas the first formally omitted error contribution has length $2n+1$.
We call this excess an \emph{overshoot}: it arises from estimating the integral remainders and is not visible in the formal order calculation.
The bound $4n-1$ is sufficient for our proof and is not claimed to be optimal: it identifies precisely which commutators must be controlled in our proof.
Once concrete operators $A$ and $B$ are specified, their commutators can be examined to select an appropriate regularity space $\mathcal{D}$, often exploiting cancellations that considerably weaken the resulting requirements.
We illustrate this mechanism for the Schr\"odinger equation with a quadratic potential.
For $A=\ii\Delta$ and $B=-\ii|x|^2$ on $L^2(\BB R^d)$, all nested commutators remain differential operators of order at most two.
Together with the regularity properties of the associated evolution groups, this ensures that the hypotheses of our theorem hold for every $n\ge1$ on the same weighted space
\begin{equation*}
	\mathcal{D}
	=\{\psi\in H^2(\BB R^d)\,\colon\,|x|^2\psi\in L^2(\BB R^d)\}.
\end{equation*}
Thus higher-order local-error estimates can be established without increasing the required regularity space.

Splittings for unbounded operators in general require to control commutators, cf.\ for example \cite{Iserles2024}.
It is perhaps germane to point out that, at least in the context of the linear Schr\"odinger equation, this is an entirely benign phenomenon.
For example, for the first-order Lie--Trotter splitting we need to control $[\Delta,V]=(\Delta V)+2(\nabla V)\cdot\nabla$ and (assuming for simplicity that $V$ is smooth) this requires $u\in H^{1+\varepsilon}$ for some $\varepsilon>0$.
However, unless $u$ has this regularity, we cannot expect the method to be of first order!
This argument extends to all the splittings in this paper, using the regularity reduction inherent in forming commutators \cite{Bader2014}.

The paper is organized as follows.
Section~\ref{sec:truncated-sbch} develops the truncated sBCH machinery in the bounded setting, including the relevant remainder representation.
In Section~\ref{sec:bounded-zassenhaus} we obtain the explicit formula for the Zassenhaus exponents, identify the leading bounded-operator local-error term, and derive the vanishing commutator identities on which our subsequent analysis depends.
Section~\ref{sec:unbounded-zassenhaus} establishes the general finite-expansion theorem for multiple conjugations of unbounded operators; the resulting commutator conditions and overshoots lead to our local-error theorem for the full symmetric Zassenhaus family and to its Schr\"odinger applications.
Additional algebraic details are collected in the appendices.

\setcounter{equation}{0}
\section{Truncated sBCH formula for bounded operators}\label{sec:truncated-sbch}
Let $A,B\colon\mathcal{H}\to\mathcal{H}$ be bounded linear operators on a given Hilbert space $\mathcal{H}$.
Let us first recall the \emph{symmetric Baker--Campbell--Hausdorff} (sBCH) series
\begin{align}\label{eq:sbch-series}
	\begin{split}
		\mathrm{sBCH}(\tau A,\tau B)
		&:=\tau A+\tau B-\frac{1}{24}\tau^{3}[[B,A],A]-\frac{1}{12}\tau^{3}[[B,A],B]\\
		&+\frac{7}{5760}\tau^{5}[[[[B,A],A],A],A]+\frac{7}{1440}\tau^{5}[[[[B,A],A],A],B]\\
		&+\frac{1}{180}\tau^{5}[[[[B,A],A],B],B]+\frac{1}{720}\tau^{5}[[[[B,A],B],B],B]\\
		&+\frac{1}{480}\tau^{5}[[[B,A],A],[B,A]]-\frac{1}{360}\tau^{5}[[[B,A],B],[B,A]]+\cdots.	
	\end{split}
\end{align} 
It is well known that
\begin{align*}
	\ee^{\frac12\tau A}\ee^{\tau B}\ee^{\frac12\tau A}
	=\ee^{\mathrm{sBCH}(\tau A,\tau B)}\quad
	\text{for sufficiently small $\tau\geq 0$;}
\end{align*}
see \cite{Blanes2024}.
Note that sBCH is an odd function of $\tau$.
For the further analysis it will be useful to introduce truncated expressions of the sBCH series.
\begin{definition}\label{def:truncated-sbch}
	For an integer $n\ge 0$, we denote by $\mathrm{sBCH}_{2n+1}(\tau A,\tau B)$ the truncation of the above infinite sBCH series up to $\mathcal{O}(\tau^{2n+3})$.
	That is, $\mathrm{sBCH}_{2n+1}(\tau A,\tau B)$ is the unique polynomial of degree at most $2n+1$ in $\tau$ coinciding with the corresponding terms of $\mathrm{sBCH}(\tau A,\tau B)$.
	Moreover, let $C_1(\tau A,\tau B):=\tau A+\tau B$,
	\begin{align*}
		C_{2n+1}(\tau A,\tau B)
		&:=\mathrm{sBCH}_{2n+1}(\tau A,\tau B)-\mathrm{sBCH}_{2n-1}(\tau A,\tau B),\quad n\ge 1,
	\end{align*}
	denote the linear combination of commutators of length $2n+1$ of $\tau A$ and $\tau B$ according to \eqref{eq:sbch-series}, whereby
	\begin{align*}
		\mathrm{sBCH}_{2n+1}^\infty(\tau A,\tau B)
		&:=\sum_{j=n}^\infty C_{2j+1}(\tau A,\tau B),\quad
		n\ge 0.
	\end{align*}
\end{definition}
Note that Definition~\ref{def:truncated-sbch} in particular yields
\begin{align}\label{eq:sbch-first-terms}
	\begin{split}
		\mathrm{sBCH}_1(\tau A,\tau B)
		&=C_1(\tau A,\tau B)
		=\tau A+\tau B,
		\quad\text{and}\\[2mm]
		\mathrm{sBCH}_3(\tau A,\tau B)
		&=\mathrm{sBCH}_1(\tau A,\tau B)+C_3(\tau A,\tau B),
		\quad\text{where}\\
		C_3(\tau A,\tau B)
		&=-\tau^{3}\frac{1}{24}[[B,A],A]-\tau^{3}\frac{1}{12}[[B,A],B].
	\end{split}
\end{align} 
\begin{remark}\label{rem:sbch-convergence}
	Note that for bounded $A$ and $B$ the convergence of the tail $\mathrm{sBCH}_{2n+1}^\infty(\tau A,\tau B)$ is guaranteed for sufficiently small $\tau>0$ by the convergence of the full sBCH series, which itself has a positive radius of convergence since it can be obtained by elementary operations of convergent series, see e.g.\ chapter III.4.2 in \cite{Hairer2006}.
	For given $A$ and $B$, we henceforth assume that $\tau$ is chosen so that the full sBCH series is convergent.
\end{remark}
In the following, we will make use of the Landau symbol.
For a one-parameter family $C(\tau)$ of bounded operators on $\mathcal{H}$ we write $C(\tau)=\mathcal{O}(\tau^k)$ for some $k\ge 0$ if there exist constants $c>0$ and $\tau_0>0$ such that, for all $\psi\in\mathcal{H}$,
\begin{align*}
	\|C(\tau)\psi\|_{\mathcal{H}}
	\le c\tau^k\|\psi\|_{\mathcal{H}}
	\quad\text{for all $0\le\tau\le\tau_0$.}
\end{align*}
In particular, for the sBCH commutators $C_{2n+1}$ we have
\begin{align}\label{eq:sbch-homogeneous-order}
	C_{2n+1}(\tau A,\tau B)
	=\mathcal{O}(\tau^{2n+1})
	\quad\text{for all $n\ge 0$.}
\end{align}
\begin{remark}
	Note that the truncated sBCH terms can also be evaluated for operators with other scalings in $\tau$.
	For instance, if $P,Q=\mathcal{O}(\tau^p)$, $p\ge 1$, we observe that
	\begin{align*}
		\mathrm{sBCH}_{2n-1}(P,Q)
		=P+Q+\mathcal{O}(\tau^{3p})
		=\mathcal{O}(\tau^p)
		\quad\text{for all $n\ge 1$},
	\end{align*}
	while for the commutators of length $2n+1$ we find
	\begin{align}\label{eq:sbch-commutator-scaling}
		C_{2n+1}(P,Q)
		=\mathcal{O}(\tau^{(2n+1)p}),
	\end{align}
	and, as a consequence, $\mathrm{sBCH}_{2n+1}^\infty(P,Q)=\mathcal{O}(\tau^{(2n+1)p})$ for all $n\ge 1$. 
\end{remark}
The following lemma identifies the leading error term caused by truncating the symmetric BCH formula.
\begin{lemma}\label{lem:sbch-truncation-error}
	Let $A,B\colon\mathcal{H}\to\mathcal{H}$ be bounded linear operators on $\mathcal{H}$.
	Then, for all $n\ge 1$, the following error representation holds:
	\begin{align*}
		\ee^{\frac12\tau A}\ee^{\tau B}\ee^{\frac12\tau A}
		=\ee^{\mathrm{sBCH}_{2n-1}(\tau A,\tau B)}+C_{2n+1}(\tau A,\tau B)+\mathcal{O}(\tau^{2n+2}).
	\end{align*}
\end{lemma}
\begin{proof}
	Since $A$ and $B$ are bounded, the sBCH series converges in operator norm for sufficiently small $\tau$; see Remark~\ref{rem:sbch-convergence}.
	Hence, for such $\tau$, the following expansion holds:
	\begin{align*}
		\ee^{\frac12\tau A}\ee^{\tau B}\ee^{\frac12\tau A}
		&=\ee^{\mathrm{sBCH}_{2n-1}(\tau A,\tau B)+\mathrm{sBCH}_{2n+1}^\infty(\tau A,\tau B)}\\
		&=\sum_{j=0}^\infty\frac{(\mathrm{sBCH}_{2n-1}(\tau A,\tau B)+\mathrm{sBCH}_{2n+1}^\infty(\tau A,\tau B))^j}{j!}\\
		&=\operatorname{Id}+\sum_{j=1}^\infty\frac{\mathrm{sBCH}_{2n-1}(\tau A,\tau B)^j}{j!}+\sum_{j=1}^\infty\frac{\mathrm{sBCH}_{2n+1}^\infty(\tau A,\tau B)^j}{j!}+\mathcal{O}(\tau^{2n+2}),
	\end{align*}
	where the scaling of the remainder $\mathcal{O}(\tau^{2n+2})$ is a consequence of the fact that it consists of products in which we have at least one factor $\mathrm{sBCH}_{2n-1}(\tau A,\tau B)=\mathcal{O}(\tau)$, and one factor $\mathrm{sBCH}_{2n+1}^\infty(\tau A,\tau B)=\mathcal{O}(\tau^{2n+1})$.
	Thus, we conclude that
	\begin{align*}
		\ee^{\frac12\tau A}\ee^{\tau B}\ee^{\frac12\tau A}
		&=\ee^{\mathrm{sBCH}_{2n-1}(\tau A,\tau B)}+\mathrm{sBCH}_{2n+1}^\infty(\tau A,\tau B)+\mathcal{O}(\tau^{2n+2})\\[3mm]
		&=\ee^{\mathrm{sBCH}_{2n-1}(\tau A,\tau B)}+C_{2n+1}(\tau A,\tau B)+\mathcal{O}(\tau^{2n+2}).
	\end{align*}
\end{proof}

Note that according to \eqref{eq:sbch-homogeneous-order} the leading-order term of the error
\begin{align}\label{eq:sbch-truncation-error}
	\begin{split}
		E_{2n-1}^{\operatorname{sBCH}}(\tau A,\tau B)
		:=&\,\ee^{\frac12\tau A}\ee^{\tau B}\ee^{\frac12\tau A}-\ee^{\mathrm{sBCH}_{2n-1}(\tau A,\tau B)}\\
		=&\,C_{2n+1}(\tau A,\tau B)+\mathcal{O}(\tau^{2n+2})	
	\end{split}
\end{align}
is $\mathcal{O}(\tau^{2n+1})$.
Moreover, for $P,Q=\mathcal{O}(\tau^p)$, $p\ge 1$, \eqref{eq:sbch-commutator-scaling} yields that
\begin{align*}
	E_{2n-1}^{\operatorname{sBCH}}(P,Q)
	=\mathcal{O}(\tau^{(2n+1)p}).	
\end{align*}
Lemma~\ref{lem:sbch-truncation-error} provides a practical error representation of the truncated sBCH formula in the case of bounded operators.
As we will show in the next section, Lemma~\ref{lem:sbch-truncation-error} serves as the point of departure for the derivation of the symmetric Zassenhaus splitting.

\setcounter{equation}{0}
\section{Symmetric Zassenhaus splittings for bounded operators}\label{sec:bounded-zassenhaus}
The symmetric Zassenhaus splitting is based on the sBCH formula introduced in the previous section and can be expressed in the form
\begin{align}\label{eq:zassenhaus-factorization}
	\ee^{\tau(A+B)}
	=\ee^{\tfrac{1}{2}\tau A}\ee^{\tfrac{1}{2}\tau B}\prod_{\ell=2}^{n}\ee^{\mathcal{R}_{2\ell-1}(\tau)}\prod_{\ell=2}^{n}\ee^{\mathcal{R}_{2(n-\ell)+3}(\tau)}\ee^{\tfrac{1}{2}\tau B}\ee^{\tfrac{1}{2}\tau A}+\mathcal{O}(\tau^{2n+1}),
\end{align}
where the \emph{Zassenhaus exponents} $\mathcal{R}_{2\ell-1}(\tau)$ are given by linear combinations of nested commutators of $A$ and $B$.
The method was first proposed in \cite{Bader2014} for the time-dependent Schr\"odinger equation in the semiclassical regime.
It was subsequently analysed in a more general framework of Banach--Lie algebras in \cite{Arnal2017}, where its convergence properties were established and a recursive construction of the operators $\mathcal{R}_{2\ell-1}$ has been derived.

In this section, we revisit the method, beginning with a presentation of the main idea underlying it, illustrated by the derivation of the first two sZS schemes with local errors $\mathcal{O}(\tau^3)$ and $\mathcal{O}(\tau^5)$.
We then provide a new formula for the Zassenhaus exponents and prove local-error bounds.

\subsection{The idea behind the general construction}\label{sec:zassenhaus-construction}

\paragraph{Local error $\mathcal{O}(\tau^3)$.}
The sZS of the lowest order coincides with the famous \emph{Strang splitting}, and can be derived by the application of the truncated sBCH formula in the following way:
\begin{align}\label{eq:strang-construction}
	\ee^{-\frac12\tau A}\ee^{\tau A+\tau B}\ee^{-\frac12\tau A}
	=\ee^{\mathrm{sBCH_1}(-\tau A,\tau A+\tau B)}+E_1^{\operatorname{sBCH}}(-\tau A,\tau A+\tau B),
\end{align}
where, according to \eqref{eq:sbch-truncation-error}, we have
\begin{align*}
	&E_1^{\operatorname{sBCH}}(-\tau A,\tau A+\tau B)
	=C_3(-\tau A,\tau A+\tau B)+\mathcal{O}(\tau^4)\\
	&\qquad=-\frac{\tau^{3}}{24}[[A+B,-A],-A]-\frac{\tau^{3}}{12}[[A+B,-A],A+B]+\mathcal{O}(\tau^4)\\
	&\qquad=\frac{\tau^{3}}{24}\Big([[B,A],A]+2[[B,A],B]\Big)+\mathcal{O}(\tau^4).
\end{align*}
Substituting $\mathrm{sBCH}_1(-\tau A,\tau A+\tau B)=\tau B$ into \eqref{eq:strang-construction} and using $\ee^{\tau X}=\operatorname{Id}+\mathcal{O}(\tau)$ for bounded $X$, we obtain
\begin{align*}
	\ee^{\tau (A+B)}
	&=\ee^{\frac12\tau A}\ee^{\tau B}\ee^{\frac12\tau A}+\ee^{\frac12\tau A}E_1^{\operatorname{sBCH}}(-\tau A,\tau A+\tau B)\ee^{\frac12\tau A}\\
	&=\ee^{\frac12\tau A}\ee^{\tau B}\ee^{\frac12\tau A}+\frac{\tau^{3}}{24}\Big([[B,A],A]+2[[B,A],B]\Big)+\mathcal{O}(\tau^4).
\end{align*}
The last equation shows that the exponential $\ee^{\tau (A+B)}$ can be approximated with local error $\mathcal{O}(\tau^3)$ by
\begin{align*}
	Z_1(\tau)
	:=\ee^{\frac12\tau A}\ee^{\tau B}\ee^{\frac12\tau A},
\end{align*}
which indeed is the well-known formula for the Strang splitting \cite{Strang1968}.

\paragraph{Local error $\mathcal{O}(\tau^5)$.}
To derive the sZS of the next higher order, we start from the second truncated sBCH formula and apply the same steps as before.
This yields 
\begin{align*}
	\ee^{-\frac12\tau A}\ee^{\tau A+\tau B}\ee^{-\frac12\tau A}
	=\ee^{\mathrm{sBCH_3}(-\tau A,\tau A+\tau B)}+E_3^{\operatorname{sBCH}}(-\tau A,\tau A+\tau B),
\end{align*}
where, according to \eqref{eq:sbch-first-terms}, we have
\begin{align*}
	\mathrm{sBCH}_3(-\tau A,\tau A+\tau B)
	=\tau B+\frac{\tau^{3}}{24}[[B,A],A]+\frac{\tau^{3}}{12}[[B,A],B].
\end{align*}
We observe that the leading-order term of $\mathrm{sBCH}_3(-\tau A,\tau A+\tau B)$ is given by $\tau B$.
Hence, splitting $\tau B$ from the remaining terms $\mathcal{O}(\tau^3)$ requires the same approach as above, \emph{i.e.},
\begin{align*}
	\ee^{-\frac12\tau B}\ee^{\mathrm{sBCH_3}(-\tau A,\tau A+\tau B)}\ee^{-\frac12\tau B}
	&=\ee^{\mathrm{sBCH}_3(-\tau B,\mathrm{sBCH}_3(-\tau A,\tau A+\tau B))}\\
	&\qquad+E_3^{\operatorname{sBCH}}(-\tau B,\mathrm{sBCH}_3(-\tau A,\tau A+\tau B)),
\end{align*}
which then leads to:
\begin{align*}
	\begin{split}
		\ee^{\tau(A+B)}
		&=\ee^{\frac12\tau A}\ee^{\frac12\tau B}\ee^{\mathrm{sBCH}_3(-\tau B,\mathrm{sBCH}_3(-\tau A,\tau A+\tau B))}\ee^{\frac12\tau B}\ee^{\frac12\tau A}\\
		&\qquad+\ee^{\frac12\tau A}E_3^{\operatorname{sBCH}}(-\tau A,\tau A+\tau B)\ee^{\frac12\tau A}\\
		&\qquad+\ee^{\frac12\tau A}\ee^{\frac12\tau B}E_3^{\operatorname{sBCH}}(-\tau B,\mathrm{sBCH}_3(-\tau A,\tau A+\tau B))\ee^{\frac12\tau B}\ee^{\frac12\tau A}.
	\end{split}
\end{align*}
In particular, a simple calculation confirms that
\begin{align*}
	\mathrm{sBCH}_3(-\tau B,\mathrm{sBCH}_3(-\tau A,\tau A+\tau B))
	=\frac{\tau^3}{24}\left([[B,A],A]+2[[B,A],B]\right)+\mathcal{O}(\tau^5),
\end{align*}
and the error terms are given by
\begin{align*}
	E_3^{\operatorname{sBCH}}(-\tau A,\tau A+\tau B)
	&=C_5(-\tau A,\tau A+\tau B)+\mathcal{O}(\tau^6),\\
	E_3^{\operatorname{sBCH}}(-\tau B,\mathrm{sBCH}_3(-\tau A,\tau A+\tau B))
	&=C_5(-\tau B,\mathrm{sBCH}_3(-\tau A,\tau A+\tau B))+\mathcal{O}(\tau^6).
\end{align*}
To conclude, we note that the above calculations show that the associated splitting
\begin{align*}
	Z_3(\tau)
	:=\ee^{\frac12\tau A}\ee^{\frac12\tau B}\ee^{\frac{\tau^3}{24}\left([[B,A],A]+2[[B,A],B]\right)}\ee^{\frac12\tau B}\ee^{\frac12\tau A}
\end{align*}
approximates the exponential $\ee^{\tau(A+B)}$ with local error $\mathcal{O}(\tau^5)$.
The construction of higher-order splittings proceeds along the same lines.
For instance, to obtain a scheme with local error $\mathcal{O}(\tau^{2n+1})$, one starts from $\mathrm{sBCH}_{2n-1}$ and then separates the leading-order contribution from the remaining higher-order terms.

Let us now take a closer look at Zassenhaus exponents.
Comparing the above expression for $Z_3$ with \eqref{eq:zassenhaus-factorization}, we conclude that the Zassenhaus exponent $\mathcal{R}_3$, which, due to symmetry, appears twice in \eqref{eq:zassenhaus-factorization}, is given by
\begin{align*}
	\mathcal{R}_3(\tau)
	=\frac{\tau^3}{48}\left([[B,A],A]+2[[B,A],B]\right).
\end{align*}
In the next step, we prove a general formula for the exponents $\mathcal{R}_{2\ell-1}$ of arbitrary order.

\subsection{A formula for Zassenhaus exponents}
Upon defining
\begin{align*}
	\mathcal{R}_{-1}(\tau)
	:=\frac12\tau A,\qquad
	\mathcal{R}_{1}(\tau)
	:=\frac12\tau B,
\end{align*}
the general form of the sZS of local error $\mathcal{O}(\tau^{2n+1})$, $n\ge 1$, presented in \eqref{eq:zassenhaus-factorization} reads
\begin{align}\label{eq:zassenhaus-splitting}
	Z_{2n-1}(\tau)
	:=\prod_{\ell=0}^{n}\ee^{\mathcal{R}_{2\ell-1}(\tau)}\prod_{\ell=0}^{n}\ee^{\mathcal{R}_{2(n-\ell)-1}(\tau)},
\end{align}
where $\mathcal{R}_{2\ell-1}(\tau)$ denotes a linear combination of commutators of length $|2\ell-1|$ of $\tau A$ and $\tau B$.
Note that, due to symmetry of the method, the term in the middle appears twice and can therefore also be written as
\begin{align*}
	\ee^{\mathcal{R}_{2n-1}(\tau)}\ee^{\mathcal{R}_{2n-1}(\tau)}
	=\ee^{2\mathcal{R}_{2n-1}(\tau)}.
\end{align*} 
As an immediate consequence, we easily observe that the elements $\mathcal{R}_{2\ell-1}(\tau)$ scale like $\tau^{|2\ell-1|}$ and that they can be derived by following the approach described in the previous subsection.
Furthermore, note that the first exponent $\mathcal{R}_{-1}(\tau)$ scales like $\tau$.

Before presenting the formula for the Zassenhaus exponents, let us introduce the following useful notation:
\begin{definition}
	For a smooth operator-valued function $F$, let
	\begin{align*}
		[s^r]F(s)
		:=\frac{1}{r!}\frac{\mathrm{d}^r}{\mathrm{d}s^r}F(s)\big|_{s=0},\qquad
		r\in\BB{Z}_{\ge0},
	\end{align*}
	denote its $r$-th Taylor coefficient at zero.
\end{definition}
We are now ready to present the formula.
\begin{theorem}\label{thm:zassenhaus-exponents}
	Let $n\ge1$.
	For a pair $\mathcal{B}=(P,Q)$, write $\mathrm{sBCH}_{2n-1}(\mathcal{B})=\mathrm{sBCH}_{2n-1}(P,Q)$ and define recursively
	\begin{align*}
		\mathcal{B}_{-1}(\tau)
		&:=(-\tau A,\tau(A+B)),\\
		\mathcal{B}_1(\tau)
		&:=\big(-\tau B,\mathrm{sBCH}_{2n-1}(\mathcal{B}_{-1}(\tau))\big),\\
		\mathcal{B}_{2\ell-1}(\tau)
		&:=\big(-2\mathcal{R}_{2\ell-1}(\tau),\mathrm{sBCH}_{2n-1}(\mathcal{B}_{2\ell-3}(\tau))\big),\qquad 2\le\ell\le n,
	\end{align*}
	where, for $2\le\ell\le n$, the Zassenhaus exponents are given by
	\begin{align}\label{eq:zassenhaus-coefficient-extraction}
		\mathcal{R}_{2\ell-1}(\tau)
		=\frac{\tau^{2\ell-1}}2[\tau^{2\ell-1}]\sum_{(j,k)\in\mathcal{J}^{2\ell-1}}C_{2j-1}(\mathcal{B}_{2k-1}(\tau)),
	\end{align}
	with $\mathcal{J}^{2\ell-1}:=\{(j,k)\in\BB{Z}^2\,\colon\,j\ge 2,\,k\ge 0,\,(2j-1)|2k-1|\le 2\ell-1\}$.\footnote{The index set can also be written as $\mathcal{J}^{2\ell-1}=\big\{(j,k)\in\BB{Z}^2\,\colon\,3\le 2j-1,\,-1\le 2k-1,\,(2j-1)\cdot|2k-1|\le 2\ell-1\big\}.$}
\end{theorem}
\begin{proof}
	We prove by induction on $\ell$ that
	\begin{align*}
		\mathrm{sBCH}_{2n-1}(\mathcal{B}_{2\ell-3}(\tau))
		=\mathcal{O}(\tau^{2\ell-1}),
	\end{align*}
	and both components of $\mathcal{B}_{2k-1}(\tau)$ are $\mathcal{O}(\tau^{|2k-1|})$ for $0\le k\le\ell-1$.
	
	The polynomials $\mathrm{sBCH}_{2n-1}(\mathcal{B}_{2k-1}(\tau))$ are odd in $\tau$.
	Initially, both components of $\mathcal{B}_{-1}$ and $\mathcal{B}_1$ are $\mathcal{O}(\tau)$, and $\mathrm{sBCH}_{2n-1}(\mathcal{B}_1(\tau))=\mathcal{O}(\tau^3)$.
	For the induction step, assume these bounds hold for some $2\le\ell\le n$.
	We introduce the index set
	\begin{align*}
		\mathcal{I}_{2n-1,2\ell-1}
		:=\{(j,k)\in\BB{Z}^2\,\colon\,2\le j\le n,\,0\le k\le\ell-1\}.
	\end{align*}
	The set $\mathcal{I}_{2n-1,2\ell-1}$ indexes all nonlinear sBCH terms accumulated before extracting $\mathcal{R}_{2\ell-1}$.
	Indeed, unfolding the recursion gives
	\begin{align*}
		\mathrm{sBCH}_{2n-1}(\mathcal{B}_{2\ell-3}(\tau))
		=-2\sum_{i=2}^{\ell-1}\mathcal{R}_{2i-1}(\tau)+\hspace{-2mm}\sum_{(j,k)\in\mathcal{I}_{2n-1,2\ell-1}}C_{2j-1}(\mathcal{B}_{2k-1}(\tau)).
	\end{align*}
	By the induction hypothesis, both components of $\mathcal{B}_{2k-1}(\tau)$ are $\mathcal{O}(\tau^{|2k-1|})$.
	Since $C_{2j-1}$ is homogeneous of degree $2j-1$ in its two arguments,
	\begin{align*}
		C_{2j-1}(\mathcal{B}_{2k-1}(\tau))
		=\mathcal{O}(\tau^{(2j-1)|2k-1|}).
	\end{align*}
	Thus the set $\mathcal{J}^{2\ell-1}$ retains the terms that may contribute to degree $2\ell-1$, whereas terms in $\mathcal{I}_{2n-1,2\ell-1}\setminus\mathcal{J}^{2\ell-1}$ start at higher degrees.
	Consequently,
	\begin{align*}
		\sum_{(j,k)\in\mathcal{I}_{2n-1,2\ell-1}\setminus\mathcal{J}^{2\ell-1}}C_{2j-1}(\mathcal{B}_{2k-1}(\tau))
		=\mathcal{O}(\tau^{2\ell+1}).
	\end{align*}
	This sum has no coefficient of degree $2\ell-1$, nor do the already extracted homogeneous exponents, since
	\begin{align*}
		[\tau^{2\ell-1}]\mathcal{R}_{2i-1}(\tau)
		=[\tau^{2\ell-1}]\big(\tau^{2i-1}\mathcal{R}_{2i-1}(1)\big)
		=0,\qquad
		2\le i\le\ell-1.
	\end{align*}
	Thus, the formula in \eqref{eq:zassenhaus-coefficient-extraction} extracts precisely the homogeneous part of the current residual:
	\begin{align*}
		2\mathcal{R}_{2\ell-1}(\tau)
		=\tau^{2\ell-1}[\tau^{2\ell-1}]\mathrm{sBCH}_{2n-1}(\mathcal{B}_{2\ell-3}(\tau)).
	\end{align*}
	By oddness, $\mathrm{sBCH}_{2n-1}(\mathcal{B}_{2\ell-3}(\tau))-2\mathcal{R}_{2\ell-1}(\tau)=\mathcal{O}(\tau^{2\ell+1})$.
	The next argument is
	\begin{align*}
		\mathcal{B}_{2\ell-1}(\tau)
		=\big(-2\mathcal{R}_{2\ell-1}(\tau),\,\mathrm{sBCH}_{2n-1}(\mathcal{B}_{2\ell-3}(\tau))\big).
	\end{align*}
	Both components of this pair are $\mathcal{O}(\tau^{2\ell-1})$, so its (nonlinear) sBCH terms start at degree $3(2\ell-1)\ge 2\ell+1$.
	Consequently, the sBCH expression satisfies
	\begin{align*}
		\mathrm{sBCH}_{2n-1}(\mathcal{B}_{2\ell-1}(\tau))
		=\mathcal{O}(\tau^{2\ell+1}),
	\end{align*}
	which completes the induction.
\end{proof}

A direct application of the formula \eqref{eq:zassenhaus-coefficient-extraction} for $Z_9$ can be found in Appendix~\ref{app:zassenhaus-exponents-z9}.
We now proceed with the analysis of the local error.

\subsection{Representation of the local error}
The next theorem establishes the error bounds for the sZS.
\begin{theorem}\label{thm:bounded-local-error}
	Let $A,B\colon\mathcal{H}\to\mathcal{H}$ be two bounded linear operators.
	For $n\ge 1$, consider the splitting $Z_{2n-1}(\tau)$ presented in \eqref{eq:zassenhaus-splitting}.
	Furthermore, define $U(\tau):=\ee^{\tau(A+B)}$.
	Then, for sufficiently small $\tau\ge0$ and all $\psi\in\mathcal{H}$, we have
	\begin{align*}
		\|U(\tau)\psi-Z_{2n-1}(\tau)\psi\|_{\mathcal{H}}
		\le (c_1\tau^{2n+1}+c_2\tau^{2n+2})\|\psi\|_\mathcal{H},
	\end{align*}
	where $c_1$ depends on commutators of length $2n+1$ of $A$ and $B$ in the following manner:
	\begin{align}\label{eq:leading-error-constant}
		c_1
		=2\|\mathcal{R}_{2n+1}(1)\|.
	\end{align}
	Moreover, $c_2>0$ depends on commutators of length $\ge 2n+2$.	
\end{theorem}
\begin{proof}
	As it was presented in Section~\ref{sec:zassenhaus-construction} for the specific cases $n=1$ and $n=2$, the first two steps in the derivation of the symmetric Zassenhaus splitting $Z_{2n-1}(\tau)$ are based on the decompositions
	\begin{align*}
		U(\tau)
		&=\ee^{\tfrac{1}{2}\tau A}\left(\ee^{\mathrm{sBCH}_{2n-1}(\mathcal{B}_{-1})}+E_{2n-1}^{\mathrm{sBCH}}(\mathcal{B}_{-1})\right)\ee^{\tfrac{1}{2}\tau A},\\
		\ee^{\mathrm{sBCH}_{2n-1}(\mathcal{B}_{-1})}
		&=\ee^{\tfrac{1}{2}\tau B}\left(\ee^{\mathrm{sBCH}_{2n-1}(\mathcal{B}_{1})}+E_{2n-1}^{\mathrm{sBCH}}(\mathcal{B}_{1})\right)\ee^{\tfrac{1}{2}\tau B}.
	\end{align*}
	Moreover, for $\ell\ge 2$, the underlying splitting reads
	\begin{align*}
		\ee^{\mathrm{sBCH}_{2n-1}(\mathcal{B}_{2\ell-3})}
		&=\ee^{\mathcal{R}_{2\ell-1}(\tau)}\left(\ee^{\mathrm{sBCH}_{2n-1}(\mathcal{B}_{2\ell-1})}+E_{2n-1}^{\mathrm{sBCH}}(\mathcal{B}_{2\ell-1})\right)\ee^{\mathcal{R}_{2\ell-1}(\tau)}.
	\end{align*}
	Let us further introduce the nested error contributions
	\begin{align*}
		\mathcal{S}_{2k-1}
		&:=\prod_{\ell=0}^{k}\ee^{\mathcal{R}_{2\ell-1}(\tau)}E_{2n-1}^{\mathrm{sBCH}}(\mathcal{B}_{2k-1})\prod_{\ell=0}^{k}\ee^{\mathcal{R}_{2(k-\ell)-1}(\tau)},\quad
		0\le k\le n.
	\end{align*}
	Note that, for ease of notation, in the following we suppress the dependence of $\mathcal{S}_{2k-1}$ on $n$.
	The combination of the individual nested splittings then produces the sZS as follows:
	\begin{align*}
		U(\tau)
		&=\ee^{\tfrac{1}{2}\tau A}\ee^{\mathrm{sBCH}_{2n-1}(\mathcal{B}_{-1})}\ee^{\tfrac{1}{2}\tau A}+\mathcal{S}_{-1}\\[1mm]
		&=\ee^{\tfrac{1}{2}\tau A}\ee^{\tfrac{1}{2}\tau B}\ee^{\mathrm{sBCH}_{2n-1}(\mathcal{B}_{1})}\ee^{\tfrac{1}{2}\tau B}\ee^{\tfrac{1}{2}\tau A}+\mathcal{S}_{-1}+\mathcal{S}_1\\[1mm]
		&\dots\\
		&=\prod_{\ell=0}^{n}\ee^{\mathcal{R}_{2\ell-1}(\tau)}\ee^{\mathrm{sBCH}_{2n-1}(\mathcal{B}_{2n-1})}\prod_{\ell=0}^{n}\ee^{\mathcal{R}_{2(n-\ell)-1}(\tau)}+\sum_{k=0}^{n}\mathcal{S}_{2k-1}.
	\end{align*}
	Applying \eqref{eq:zassenhaus-coefficient-extraction} with $\ell=n+1$ requires the recursion using $\mathrm{sBCH}_{2n+1}$, whereas the arguments $\mathcal{B}_{2k-1}$ above were constructed using $\mathrm{sBCH}_{2n-1}$.
	Since
	\begin{align*}
		\mathrm{sBCH}_{2n+1}(P,Q)
		=\mathrm{sBCH}_{2n-1}(P,Q)+C_{2n+1}(P,Q),
	\end{align*}
	both recursions produce the same exponents up to and including $\mathcal{R}_{2n-1}$, but the original recursion omits contributions of degree $2n+1$ in $\tau$.
	However, only the terms $C_{2n+1}(\mathcal{B}_{-1})$ and $C_{2n+1}(\mathcal{B}_1)$ contribute at this degree.
	Later arguments are $\mathcal{O}(\tau^3)$, and changes of order $\tau^{2n+1}$ in existing nonlinear terms contribute only at degree $2n+3$ or higher.
	Restoring these two contributions, \eqref{eq:zassenhaus-coefficient-extraction} with $\ell=n+1$ therefore gives
	\begin{align*}
		2\mathcal{R}_{2n+1}(\tau)
		=\tau^{2n+1}[\tau^{2n+1}]\Big(\mathrm{sBCH}_{2n-1}(\mathcal{B}_{2n-1}(\tau))+C_{2n+1}(\mathcal{B}_{-1}(\tau))+C_{2n+1}(\mathcal{B}_1(\tau))\Big).
	\end{align*}
	Oddness and the order of the final sBCH polynomial then yield
	\begin{align*}
		\mathrm{sBCH}_{2n-1}(\mathcal{B}_{2n-1})
		=2\mathcal{R}_{2n+1}(\tau)-\tau^{2n+1}[\tau^{2n+1}]\Big(C_{2n+1}(\mathcal{B}_{-1}(\tau))+C_{2n+1}(\mathcal{B}_1(\tau))\Big)+\mathcal{O}(\tau^{2n+3}).
	\end{align*}
	Using the fact that each of the exponentials in
	\begin{align*}
		\prod_{\ell=0}^{n}\ee^{\mathcal{R}_{2\ell-1}(\tau)}\ee^{\mathrm{sBCH}_{2n-1}(\mathcal{B}_{2n-1})}\prod_{\ell=0}^{n}\ee^{\mathcal{R}_{2(n-\ell)-1}(\tau)}
	\end{align*}
	can be expanded as $\ee^{X}=\operatorname{Id}+X+\mathcal{O}(X^2)$, this results in
	\begin{align*}
		U(\tau)
		&=Z_{2n-1}(\tau)+2\mathcal{R}_{2n+1}(\tau)+\sum_{k=0}^{n}\mathcal{S}_{2k-1}\\
		&\qquad-\tau^{2n+1}[\tau^{2n+1}]\Big(C_{2n+1}(\mathcal{B}_{-1}(\tau))
		+C_{2n+1}(\mathcal{B}_1(\tau))\Big)+\mathcal{O}(\tau^{2n+2}).
	\end{align*}
	Furthermore, combining \eqref{eq:sbch-truncation-error} with the same expansion of the exponentials appearing in $\mathcal{S}_{2k-1}$ yields
	\begin{align*}
		\mathcal{S}_{2k-1}
		&=E_{2n-1}^{\mathrm{sBCH}}(\mathcal{B}_{2k-1})+\mathcal{O}(\tau^{2n+2})
		=C_{2n+1}(\mathcal{B}_{2k-1})+\mathcal{O}(\tau^{2n+2}),\quad
		k\in\{0,1\},
	\end{align*}
	as well as $\mathcal{S}_{2k-1}=\mathcal{O}(\tau^{(2n+1)(2k-1)})$ for all $k\ge 2$.
	The contributions of $\mathcal{S}_{-1}$ and $\mathcal{S}_1$ cancel the subtracted coefficient in the preceding expression for $U(\tau)$.
	Thus
	\begin{align*}
		U(\tau)-Z_{2n-1}(\tau)
		&=2\mathcal{R}_{2n+1}(\tau)+\mathcal{O}(\tau^{2n+2}).
	\end{align*}
	Since $\mathcal{R}_{2n+1}(\tau)=\tau^{2n+1}\mathcal{R}_{2n+1}(1)$, the leading coefficient is bounded by $c_1$ in \eqref{eq:leading-error-constant}.
	Note that the total local error $E_{2n-1}(\tau)$ is in fact given by an infinite series in powers of $\tau^{2n+k}$, $k\ge 1$, and convergence is guaranteed by the general assumption that $\tau$ was chosen such that the full sBCH series converges, see Remark~\ref{rem:sbch-convergence}.
	Consequently, for all $\psi\in\mathcal{H}$ we conclude that
	\begin{align*}
		\|E_{2n-1}(\tau)\psi\|_{\mathcal{H}}
		\le (c_1\tau^{2n+1}+c_2\tau^{2n+2})\|\psi\|_\mathcal{H},
	\end{align*}
	with $c_1$ in \eqref{eq:leading-error-constant}, and $c_2>0$ depending on commutators of length larger than $2n+1$ of $A$ and $B$.
\end{proof}

The smoothness of the error function for bounded $A$ and $B$ yields the following result.
\begin{corollary}\label{cor:vanishing-error-derivatives}
	Under the assumptions of Theorem~\ref{thm:bounded-local-error}, the local-error function
	\begin{align}\label{eq:local-error-definition}
		E_{2n-1}(\tau)
		=U(\tau)-Z_{2n-1}(\tau)
	\end{align}
	is smooth with vanishing derivatives at zero up to order $2n$, \emph{i.e.},
	\begin{align*}
		\frac{\mathrm{d}^s}{\mathrm{d}\tau^s}E_{2n-1}(\tau)\big|_{\tau=0}
		=0,\quad
		s
		=0,1,\dots,2n.
	\end{align*}
\end{corollary}
\begin{proof}
	Since $A$, $B$, and the Zassenhaus exponents are bounded operators, both $U(\tau)=\ee^{\tau(A+B)}$ and $Z_{2n-1}(\tau)$ are analytic functions of $\tau$.
	Consequently, $E_{2n-1}=U-Z_{2n-1}$ is analytic.
	Moreover, Theorem~\ref{thm:bounded-local-error} yields $\|E_{2n-1}(\tau)\|=\mathcal{O}(\tau^{2n+1})$.
	The Taylor expansion of $E_{2n-1}$ at zero therefore has vanishing coefficients up to and including degree $2n$, and hence $E_{2n-1}^{(j)}(0)=0$, for all $0\le j\le 2n$.
\end{proof}
%

\subsection{Integral representation of the local error}
In the previous subsection we presented the local error for the sZS in the case of bounded operators.
In this subsection we introduce an alternative integral representation of the aforementioned local error based on the Duhamel formula.
Together with Corollary~\ref{cor:vanishing-error-derivatives}, this formulation provides the basis for our main result: a local-error bound of order $\mathcal{O}(\tau^{2n+1})$ for the symmetric Zassenhaus splitting $Z_{2n-1}$ with unbounded operators, under explicit regularity assumptions.

For ease of notation, from now on we suppress the dependence of $Z_{2n-1}^{[j,k]}(\tau)$ on $2n-1$.
We start with the following definition:
\begin{definition}
	For $n\ge 1$ and odd integers $j,k$ with $-1\le j\le k\le 2n-1$ we denote by $Z^{[j,k]}(\tau)$ the product of subsequent Zassenhaus factors in the sZS in \eqref{eq:zassenhaus-factorization}, starting from the $j$-th and ending at the $k$-th, \emph{i.e.},
	\begin{align*}
		Z^{[j,k]}(\tau)
		:=\ee^{\mathcal{R}_j(\tau)}\ee^{\mathcal{R}_{j+2}(\tau)}\ldots\ee^{\mathcal{R}_{k-2}(\tau)}\ee^{\mathcal{R}_k(\tau)}.
	\end{align*}
	Moreover, we introduce the flipped-index notation
	\begin{align*}
		Z^{[k,j]}(\tau)
		:=\ee^{\mathcal{R}_k(\tau)}\ee^{\mathcal{R}_{k-2}(\tau)}\ldots\ee^{\mathcal{R}_{j+2}(\tau)}\ee^{\mathcal{R}_j(\tau)}=\left(Z^{[j,k]}(-\tau)\right)^{-1}\hspace{-2mm}.
	\end{align*}
\end{definition}
Note that the previous definition in particular yields the following short-hand expressions:
\begin{align*}
	Z^{[-1,-1]}(\tau)
	&=\ee^{\tfrac{1}{2}\tau A},\quad
	Z^{[-1,1]}(\tau)
	=\ee^{\tfrac{1}{2}\tau A}\ee^{\tfrac{1}{2}\tau B},\quad
	Z^{[1,-1]}(\tau)
	=\ee^{\tfrac{1}{2}\tau B}\ee^{\tfrac{1}{2}\tau A},\\[-8mm]
\end{align*}
\begin{align*}
	Z^{[j,k]}(\tau)Z^{[k,j]}(-\tau)
	&=\operatorname{Id},\quad
	Z_{2n-1}(\tau)
	=Z^{[-1,2n-1]}(\tau)Z^{[2n-1,-1]}(\tau).
\end{align*}
Following the approach of \cite{Lasser2020}, we derive a two-sided variation-of-constants representation of the local error, extending their construction to arbitrary-order symmetric Zassenhaus splittings.
\begin{lemma}\label{lem:local-error-representation}
	The local-error function in \eqref{eq:local-error-definition} satisfies the differential equation
	\begin{align}\label{eq:local-error-equation}
		\frac{\mathrm{d}}{\mathrm{d}\tau}E_{2n-1}(\tau)
		=\frac{(A+B)}{2}E_{2n-1}(\tau)+E_{2n-1}(\tau)\frac{(A+B)}{2}+N_{2n-1}(\tau),
	\end{align}
	with the initial condition $E_{2n-1}(0)=0$ and the remainder
	\begin{align}\label{eq:defect-definition}
		\begin{split}
			N_{2n-1}(\tau)
			&:=\frac12BZ_{2n-1}(\tau)+\frac12Z_{2n-1}(\tau)B\\
			&\qquad-\sum_{\ell=1}^{n}Z^{[-1,2\ell-3]}(\tau)\frac{\mathrm{d}}{\mathrm{d}\tau}\mathcal{R}_{2\ell-1}(\tau)Z^{[2\ell-3,-1]}(-\tau)Z_{2n-1}(\tau)\\
			&\qquad-\sum_{\ell=1}^{n}Z_{2n-1}(\tau)Z^{[-1,2\ell-3]}(-\tau)\frac{\mathrm{d}}{\mathrm{d}\tau}\mathcal{R}_{2\ell-1}(\tau)Z^{[2\ell-3,-1]}(\tau).
		\end{split}
	\end{align}
	Moreover, the error can be written in the compact integral form
	\begin{align}\label{eq:duhamel-error-representation}
		E_{2n-1}(\tau)
		=\int_0^\tau U\left(\tfrac{1}2(\tau-s)\right)N_{2n-1}(s)U\left(\tfrac{1}2(\tau-s)\right)\,\mathrm{d}s.
	\end{align}
\end{lemma}
\begin{proof}
	The derivative of the error function is computed using the product rule as
	\begin{align*}
		\frac{\mathrm{d}}{\mathrm{d}\tau}E_{2n-1}(\tau)
		&=\frac{\mathrm{d}}{\mathrm{d}\tau}U(\tau)-\frac{\mathrm{d}}{\mathrm{d}\tau}Z_{2n-1}(\tau)
		=\frac{(A+B)}{2}U(\tau)+U(\tau)\frac{(A+B)}{2}\\
		&\qquad-\sum_{\ell=0}^{n}Z^{[-1,2\ell-3]}(\tau)\frac{\mathrm{d}}{\mathrm{d}\tau}\mathcal{R}_{2\ell-1}(\tau)Z^{[2\ell-3,-1]}(-\tau)Z_{2n-1}(\tau)\\
		&\qquad-\sum_{\ell=0}^{n}Z_{2n-1}(\tau)Z^{[-1,2\ell-3]}(-\tau)\frac{\mathrm{d}}{\mathrm{d}\tau}\mathcal{R}_{2\ell-1}(\tau)Z^{[2\ell-3,-1]}(\tau),
	\end{align*}
	where, for $\ell=0$, we introduced the definition $Z^{[-1,-3]}(\pm\tau):=\operatorname{Id}$.
	Moreover, we have
	\begin{align*}
		\frac{(A+B)}{2}U(\tau)+U(\tau)\frac{(A+B)}{2}
		&=\frac{(A+B)}{2}E_{2n-1}(\tau)+E_{2n-1}(\tau)\frac{(A+B)}{2}\\
		&\qquad+\frac{(A+B)}{2}Z_{2n-1}(\tau)+Z_{2n-1}(\tau)\frac{(A+B)}{2}.
	\end{align*}
	We therefore obtain
	\begin{align*}
		\frac{\mathrm{d}}{\mathrm{d}\tau}E_{2n-1}(\tau)
		&=\frac{(A+B)}{2}E_{2n-1}(\tau)+E_{2n-1}(\tau)\frac{(A+B)}{2}\\
		&\qquad+\frac12BZ_{2n-1}(\tau)+\frac12Z_{2n-1}(\tau)B\\
		&\qquad-\sum_{\ell=1}^{n}Z^{[-1,2\ell-3]}(\tau)\frac{\mathrm{d}}{\mathrm{d}\tau}\mathcal{R}_{2\ell-1}(\tau)Z^{[2\ell-3,-1]}(-\tau)Z_{2n-1}(\tau)\\
		&\qquad-\sum_{\ell=1}^{n}Z_{2n-1}(\tau)Z^{[-1,2\ell-3]}(-\tau)\frac{\mathrm{d}}{\mathrm{d}\tau}\mathcal{R}_{2\ell-1}(\tau)Z^{[2\ell-3,-1]}(\tau)\\[1mm]
		&=\frac{(A+B)}{2}E_{2n-1}(\tau)+E_{2n-1}(\tau)\frac{(A+B)}{2}+N_{2n-1}(\tau),	
	\end{align*}
	which shows \eqref{eq:local-error-equation}.
	It therefore suffices to prove \eqref{eq:duhamel-error-representation}. 
	Indeed, direct differentiation yields
	\begin{align*}
		&\frac{\mathrm{d}}{\mathrm{d}\tau}\int_0^\tau U\left(\tfrac{1}2(\tau-s)\right)N_{2n-1}(s)U\left(\tfrac{1}2(\tau-s)\right)\,\mathrm{d}s\\
		&\qquad=U(0)N_{2n-1}(\tau)U(0)+\int_0^\tau\frac{\partial}{\partial\tau}\left[U\left(\tfrac{1}2(\tau-s)\right)N_{2n-1}(s)U\left(\tfrac{1}2(\tau-s)\right)\right]\mathrm{d}s\\
		&\qquad=N_{2n-1}(\tau)+\int_0^\tau\frac{(A+B)}{2}U\left(\tfrac{1}2(\tau-s)\right)N_{2n-1}(s)U\left(\tfrac{1}2(\tau-s)\right)\,\mathrm{d}s\\
		&\qquad\qquad+\int_0^\tau U\left(\tfrac{1}2(\tau-s)\right)N_{2n-1}(s)U\left(\tfrac{1}2(\tau-s)\right)\frac{(A+B)}{2}\,\mathrm{d}s.
	\end{align*}
The integral satisfies \eqref{eq:local-error-equation} and vanishes at $\tau=0$.
Uniqueness for this bounded linear differential equation therefore proves \eqref{eq:duhamel-error-representation}.
\end{proof}
%

\subsection{Simplification of algebraic commutator identities}
Recall that in Lemma~\ref{lem:local-error-representation} we derived an alternative integral representation of the local-error function of the splitting $Z_{2n-1}$, expressed in terms of the remainder function $N_{2n-1}$ introduced in \eqref{eq:defect-definition}.
By Corollary~\ref{cor:vanishing-error-derivatives}, the first $2n$ derivatives of the local-error function $E_{2n-1}$ vanish identically at $\tau=0$.
Consequently, the derivatives of order $0,\dots,2n-1$ of the remainder function $N_{2n-1}$ must also vanish at zero.
As shown below, these cancellations yield simplifications of algebraic identities: the Lie polynomials arising as coefficients of $\tau^k$, $0\le k\le2n-1$, in a suitable factorization of the remainder vanish identically, thereby ensuring $N_{2n-1}(\tau)=\mathcal{O}(\tau^{2n})$.
These identities will also be used in the unbounded setting.

Let us proceed by first presenting an alternative representation of the remainder function.
\begin{lemma}\label{lem:defect-decomposition}
	The remainder function $N_{2n-1}$ defined in \eqref{eq:defect-definition} can be written as
	\begin{align}\label{eq:defect-decomposition}
		N_{2n-1}(\tau)
		=\sum_{\ell=1}^n\Big[Z_{2n-1}(\tau)F_\ell(\tau)+F_\ell(-\tau)Z_{2n-1}(\tau)\Big]\tau^{2\ell-2},
	\end{align}
	with functions $F_\ell$ given by
	\begin{align*}
		F_1(\tau)
		&:=\frac12B-\frac12\ee^{-\tfrac{1}{2}\tau A}B\ee^{\tfrac{1}{2}\tau A},\\
		F_\ell(\tau)
		&:=-(2\ell-1)Z_{2n-1}^{[-1,2\ell-3]}(-\tau)\mathcal{R}_{2\ell-1}(1)Z_{2n-1}^{[2\ell-3,-1]}(\tau),\quad
		2\le\ell\le n.
	\end{align*}
\end{lemma}
\begin{proof}
	For all $\ell\ge 1$ we have $\mathcal{R}_{2\ell-1}(\tau)=\tau^{2\ell-1}\mathcal{R}_{2\ell-1}(1)$ and conclude that
	\begin{align*}
		\frac{\mathrm{d}}{\mathrm{d}\tau}\mathcal{R}_{2\ell-1}(\tau)
		=(2\ell-1)\tau^{2\ell-2}\mathcal{R}_{2\ell-1}(1).
	\end{align*}
	Hence, \eqref{eq:defect-definition} takes the form
	\begin{align*}
		N_{2n-1}(\tau)
		=&\,\,\frac12BZ_{2n-1}(\tau)+\frac12Z_{2n-1}(\tau)B\\
		&\,-\frac12\ee^{\tfrac{1}{2}\tau A}B\ee^{-\tfrac{1}{2}\tau A}Z_{2n-1}(\tau)-\frac12Z_{2n-1}(\tau)\ee^{-\tfrac{1}{2}\tau A}B\ee^{\tfrac{1}{2}\tau A}\\
		&\,-\sum_{\ell=2}^{n}Z^{[-1,2\ell-3]}(\tau)\mathcal{R}_{2\ell-1}(1)Z^{[2\ell-3,-1]}(-\tau)Z_{2n-1}(\tau)(2\ell-1)\tau^{2\ell-2}\\
		&\,-\sum_{\ell=2}^{n}Z_{2n-1}(\tau)Z^{[-1,2\ell-3]}(-\tau)\mathcal{R}_{2\ell-1}(1)Z^{[2\ell-3,-1]}(\tau)(2\ell-1)\tau^{2\ell-2}
	\end{align*}
	and a comparison with \eqref{eq:defect-decomposition} yields the stated form of the functions $F_\ell$. 
\end{proof}

In the next step we will derive expansions of the remainder function $N_{2n-1}$ based on expansions of conjugated operators of the form $\ee^{\tau^kX}Y\ee^{-\tau^kX}$.
As we can see, the functions $F_\ell$ in \eqref{eq:defect-decomposition} are linear combinations of nested conjugations of the more general form
\begin{align}\label{eq:nested-conjugation-form}
	\ee^{\tau^{k_n}X_n}\cdots\ee^{\tau^{k_1}X_1}Y\ee^{-\tau^{k_1}X_1}\cdots\ee^{-\tau^{k_n}X_n}.
\end{align}
Such expressions can be expanded by applying the single-conjugation expansion repeatedly.
To preserve this structure in the remainder, we factor out $Z_{2n-1}(\tau)$:
\begin{align}\label{eq:factored-defect}
	N_{2n-1}(\tau)
	&=G_{2n-1}(\tau)Z_{2n-1}(\tau),\\
	G_{2n-1}(\tau)
	&:=\sum_{\ell=1}^n\tau^{2\ell-2}\Big(Z_{2n-1}(\tau)F_\ell(\tau)Z_{2n-1}(-\tau)+F_\ell(-\tau)\Big).\label{eq:factored-defect-coefficients}
\end{align}
Indeed, both terms $Z_{2n-1}(\tau)F_\ell(\tau)Z_{2n-1}(-\tau)$ and $F_\ell(-\tau)$ are linear combinations of nested conjugations of the form \eqref{eq:nested-conjugation-form}.
We therefore expand $G_{2n-1}$ and recover the remainder $N_{2n-1}$ by multiplication with $Z_{2n-1}$.
Before deriving these expansions, we simplify the conjugations $Z_{2n-1}(\tau)F_\ell(\tau)Z_{2n-1}(-\tau)$ by cancelling adjacent exponential factors.

For ease of notation, we suppress the subscripts $2n-1$.
The existing product notation makes the cancellations explicit:
\begin{align*}
	Z(\tau)Z^{[-1,2\ell-3]}(-\tau)
	&=Z^{[-1,2n-1]}(\tau)Z^{[2n-1,2\ell-1]}(\tau),\\
	Z^{[2\ell-3,-1]}(\tau)Z(-\tau)
	&=Z^{[2\ell-1,2n-1]}(-\tau)Z^{[2n-1,-1]}(-\tau).
\end{align*}
The factors adjacent to $\mathcal{R}_{2\ell-1}(1)$ then cancel, since
\begin{align*}
	\ee^{\mathcal{R}_{2\ell-1}(\tau)}\mathcal{R}_{2\ell-1}(1)\ee^{-\mathcal{R}_{2\ell-1}(\tau)}
	=\mathcal{R}_{2\ell-1}(1).
\end{align*}
Consequently, for $2\le\ell\le n$,
\begin{align}\label{eq:simplified-factored-defect}
	Z(\tau)F_\ell(\tau)Z(-\tau)
	&=-(2\ell-1)Z^{[-1,2n-1]}(\tau)Z^{[2n-1,2\ell+1]}(\tau)\mathcal{R}_{2\ell-1}(1)\nonumber\\
	&\qquad\times Z^{[2\ell+1,2n-1]}(-\tau)Z^{[2n-1,-1]}(-\tau).
\end{align}
Here $Z^{[2n-1,2n+1]}(\tau)=Z^{[2n+1,2n-1]}(-\tau):=\operatorname{Id}$ denote empty products.
For $F_1$, the same cancellation uses $\mathcal{R}_1(\tau)=\tau B/2$ and gives
\begin{align*}
	Z(\tau)F_1(\tau)Z(-\tau)
	&=\frac12 Z(\tau)BZ(-\tau)-\frac12 Z^{[-1,2n-1]}(\tau)Z^{[2n-1,3]}(\tau)B\\
	&\qquad\times Z^{[3,2n-1]}(-\tau)Z^{[2n-1,-1]}(-\tau).
\end{align*}
For $\ell=n\ge2$, the remaining rightmost factor also commutes with $\mathcal{R}_{2n-1}(1)$, so
\begin{align*}
	Z(\tau)F_n(\tau)Z(-\tau)
	=-(2n-1)Z^{[-1,2n-3]}(\tau)\mathcal{R}_{2n-1}(1)Z^{[2n-3,-1]}(-\tau).
\end{align*}
Thus the highest exponent disappears from this contribution.

The above formulas express the factored remainder $G$ as a sum of nested conjugations.
We next derive their expansions, beginning with a single conjugation and then treating successive conjugations.
First, we introduce the notation for iterated commutators.
\begin{definition}
	For linear operators $X$ and $Y$, not necessarily bounded, we recall the definition of the iterated commutators
	\begin{align*}
		\operatorname{ad}_X^0(Y)
		:=Y,\quad
		\operatorname{ad}_X^i(Y)
		:=[X,\operatorname{ad}_X^{i-1}(Y)],\quad
		\text{for all $i\ge 1$}.
	\end{align*}
\end{definition}
In the presence of bounded operators $X$ and $Y$, we also recall the well-known identity (see e.g.\ Proposition~3.35 in \cite{Hall2015})
\begin{align*}
	\ee^{tX}Y\ee^{-tX}
	=\ee^{t\operatorname{ad}_X}Y
	=\sum_{j=0}^\infty\frac{t^j}{j!}\operatorname{ad}^j_X(Y),
\end{align*}
which converges in operator norm for every $t\in\BB{R}$.

We are now ready to establish the structure of the power series expansion of nested time-reversible products as they appear in the factored remainder $G$ from \eqref{eq:factored-defect-coefficients}.
In Proposition~\ref{prop:bounded-conjugation-expansion} and Theorem~\ref{thm:multiple-conjugation-expansion} below, $n$ denotes the number of conjugations rather than the splitting order, and the symbols $A_n(t)$, $B_{k,n}$, and $B_{p,n}$ are unrelated to the operators $A$ and $B$.
\begin{proposition}\label{prop:bounded-conjugation-expansion}
	Let $n\ge 1$.
	For bounded operators $A$ and $B$, we assume that each $X_j$, $1\le j\le n$, and $Y$ is a homogeneous Lie polynomial in $A,B$.
	Furthermore, let $\ell(\cdot)$ denote the length of a given commutator, assume that $\ell(X_j)=k_j$, and let
	\begin{align*}
		A_n(t)
		:=\left(\ee^{t^{k_n}X_n}\cdots\ee^{t^{k_1}X_1}\right)Y\left(\ee^{-t^{k_1}X_1}\cdots\ee^{-t^{k_n}X_n}\right).
	\end{align*}
	Then there exist operators $B_{k,n}$, each either identically zero or a linear combination of nested commutators of length $k+\ell(Y)$, such that, for every $t\in\BB{R}$, we have
	\begin{align}\label{eq:bounded-conjugation-expansion}
		A_n(t)
		=\sum_{k=0}^\infty B_{k,n}t^k.
	\end{align}
\end{proposition}
\begin{proof}
	Let us observe that, upon defining $A_0(t):=Y$, we have
	\begin{align*}
		A_n(t)
		=\sum_{j_n=0}^\infty\frac{t^{j_nk_n}}{j_n!}\operatorname{ad}^{j_n}_{X_n}\left(A_{n-1}(t)\right),\quad
		n\ge 1.
	\end{align*}	
	Iterating this expansion $n$ times and using the linearity of the adjoint operator yields
	\begin{align*}
		A_n(t)
		&=\sum_{j_n=0}^\infty\frac{t^{j_nk_n}}{j_n!}\operatorname{ad}^{j_n}_{X_n}\left(\sum_{j_{n-1}=0}^\infty\frac{t^{j_{n-1}k_{n-1}}}{j_{n-1}!}\operatorname{ad}^{j_{n-1}}_{X_{n-1}}\left(\dots\sum_{j_1=0}^\infty\frac{t^{j_1k_1}}{j_1!}\operatorname{ad}^{j_1}_{X_1}\left(Y\right)\right)\right)\\
		&=\sum_{j_n=0}^\infty\sum_{j_{n-1}=0}^\infty\cdots\sum_{j_1=0}^\infty\frac{t^{j_nk_n+j_{n-1}k_{n-1}+\dots+j_1k_1}}{j_1!j_{2}!\cdots j_n!}\operatorname{ad}^{j_n}_{X_n}\operatorname{ad}^{j_{n-1}}_{X_{n-1}}\cdots\operatorname{ad}^{j_1}_{X_1}(Y)\\
		&=\sum_{k=0}^\infty\left(\sum_{j\in\Omega_{k,n}}\frac{1}{j!}\operatorname{ad}^{j_n}_{X_n}\operatorname{ad}^{j_{n-1}}_{X_{n-1}}\cdots\operatorname{ad}^{j_1}_{X_1}(Y)\right)t^k,
	\end{align*}
	where we introduced the multi-index set
	\begin{align*}
		\Omega_{k,n}
		:=\big\{j=(j_1,\dots,j_n)\in\BB{Z}_{\ge 0}^n\,\colon\,j_1k_1+\dots+j_nk_n=k\big\}
	\end{align*}	
	and $j!=j_1!\cdots j_n!$ for the corresponding multi-index.
	The coefficients are the linear combinations of nested commutators
	\begin{align*}
		B_{k,n}
		=\sum_{j\in\Omega_{k,n}}\frac{1}{j!}\operatorname{ad}^{j_n}_{X_n}\operatorname{ad}^{j_{n-1}}_{X_{n-1}}\cdots\operatorname{ad}^{j_1}_{X_1}(Y).
	\end{align*}
	Each summand is homogeneous of length $j_n\ell(X_n)+\dots+j_1\ell(X_1)+\ell(Y)=k+\ell(Y)$ for $j\in\Omega_{k,n}$, so every nonzero $B_{k,n}$ has this same length.
\end{proof}
\begin{remark}\label{rem:defect-component-expansions}
	According to Proposition~\ref{prop:bounded-conjugation-expansion}, the functions $F_\ell$ introduced in Lemma~\ref{lem:defect-decomposition} have the expansion
	\begin{align*}
		F_\ell(\tau)
		=\sum_{k=0}^\infty B_{k,\ell}\tau^k,\quad
		1\le\ell\le n,
	\end{align*}
	where the operators $B_{k,\ell}$ are either identically zero or linear combinations of nested commutators of $A,B$, all of length $k+2\ell-1$.
	In particular, the first two functions possess the expansions
	\begin{align*}
		F_1(\tau)
		&=\frac12B-\frac12\ee^{-\tfrac{1}{2}\tau A}B\ee^{\tfrac{1}{2}\tau A}\\
		&=\frac12B-\frac12\left(B-\frac{\tau}{2}[A,B]+\frac{\tau^2}{8}[A,[A,B]]+\mathcal{O}(\tau^3)\right)\\
		&=\frac{1}{4}[A,B]\,\tau-\frac{1}{16}[A,[A,B]]\,\tau^2+\mathcal{O}(\tau^3),\\[2mm]
		F_2(\tau)
		&=-3Z^{[-1,1]}(-\tau)\mathcal{R}_{3}(1)Z^{[1,-1]}(\tau)\\[1mm]
		&=-\frac{3}{48}\ee^{-\frac12\tau A}\ee^{-\frac12\tau B}\left([[B,A],A]+2[[B,A],B]\right)\ee^{\frac12\tau B}\ee^{\frac12\tau A}\\
		&=-\frac{3}{48}\left([[B,A],A]+2[[B,A],B]\right)+\mathcal{O}(\tau),
	\end{align*}
	which demonstrates that
	\begin{align*}
		B_{0,1}
		=0,\quad
		B_{1,1}
		=\frac{1}{4}[A,B],\quad
		B_{2,1}
		=-\frac{1}{16}[A,[A,B]],
	\end{align*}
	as well as
	\begin{align*}
		B_{0,2}
		=-\frac{1}{16}\left([A,[A,B]]-2[[A,B],B]\right).
	\end{align*}
\end{remark}
We now expand the factored defect $G$, so that its coefficients are fixed Lie polynomials; recall that the subscripts $2n-1$ are suppressed.
\begin{lemma}\label{lem:factored-defect-coefficients}
	For bounded linear $A,B$, the factored defect $G$ defined in \eqref{eq:factored-defect-coefficients} has the convergent expansion
	\begin{align}\label{eq:factored-defect-series}
		G(\tau)
		=\sum_{k=0}^{\infty}K_k\tau^k,\qquad
		K_k
		:=[\tau^k]G(\tau).
	\end{align}
	In particular, each coefficient $K_k$ is zero or a homogeneous Lie polynomial of length $k+1$ in $A,B$.
\end{lemma}
\begin{proof}
	{
	By \eqref{eq:simplified-factored-defect} and the corresponding formula for $F_1$, the $\ell$-th summand in $G$ is $\tau^{2\ell-2}$ times a linear combination of conjugations of $\mathcal{R}_{2\ell-1}(1)$ for $\ell\ge2$, or of $B$ for $\ell=1$.
	In both cases, the conjugated polynomial is homogeneous of commutator length $2\ell-1$.
	Proposition~\ref{prop:bounded-conjugation-expansion} expands each conjugation as a power series in $\tau$.
	Let $p\ge 0$ denote the power of $\tau$ in a term of this expansion.
	Its coefficient is a homogeneous Lie polynomial of length $p+2\ell-1$.
	Multiplication by $\tau^{2\ell-2}$ changes the power to $k=p+2\ell-2$ without changing the coefficient.
	Consequently, this coefficient has length
	\begin{align*}
		p+2\ell-1
		=k+1.	
	\end{align*}
	Thus every coefficient $K_k$ is either zero or a homogeneous Lie polynomial of length $k+1$.
	}
\end{proof}

Corollary~\ref{cor:vanishing-error-derivatives} provides the cancellations needed for the following simplification of algebraic commutator identities, which will also be used in the unbounded setting.
\begin{theorem}\label{thm:vanishing-defect-coefficients}
	The coefficients in \eqref{eq:factored-defect-series} satisfy
	\begin{align}\label{eq:vanishing-defect-coefficients}
		K_k
		=[\tau^k]G(\tau)
		=0,\qquad
		0\le k\le 2n-1.
	\end{align}
\end{theorem}
\begin{proof}
	Write $E(\tau)=U(\tau)-Z(\tau)$.
	By Corollary~\ref{cor:vanishing-error-derivatives}, $E$ is smooth and $E^{(j)}(0)=0$ for $0\le j\le 2n$.
	Taylor's theorem therefore gives $E(\tau)=\mathcal{O}(\tau^{2n+1})$, $\tfrac{\mathrm{d}}{\mathrm{d}\tau}E(\tau)=\mathcal{O}(\tau^{2n})$.
	The defect equation \eqref{eq:local-error-equation} yields
	\begin{align*}
		N(\tau)
		=\frac{\mathrm{d}}{\mathrm{d}\tau}E(\tau)-\frac{(A+B)}{2}E(\tau)-E(\tau)\frac{(A+B)}{2}
		=\mathcal{O}(\tau^{2n}).
	\end{align*}
	Since $Z(-\tau)$ is locally bounded, \eqref{eq:factored-defect} gives
	\begin{align*}
		G(\tau)
		=N(\tau)Z(-\tau)
		=\mathcal{O}(\tau^{2n}).
	\end{align*}
	Comparing Taylor coefficients proves \eqref{eq:vanishing-defect-coefficients}.
\end{proof}

We note that the cancellations in Theorem~\ref{thm:vanishing-defect-coefficients} are algebraic: after expanding the commutators and collecting identical products of $A$ and $B$, every coefficient vanishes.
Thus the simplification of algebraic commutator identities holds independently of the particular operators, wherever the required products are defined.
\begin{example}
	We establish the first three vanishing identities explicitly:
	
	\noindent\textbf{First identity ($k=0$):}
	Only $\ell=1$ contributes to the constant coefficient in \eqref{eq:factored-defect-coefficients}, since $\tau^{2\ell-2}$ has positive degree for $\ell\ge2$.
	Using $F_1(0)=B_{0,1}=0$ and $Z(0)=\operatorname{Id}$ gives
	\begin{align*}
		K_0
		=G(0)
		=\operatorname{Id}\,B_{0,1}\operatorname{Id}+B_{0,1}
		=2B_{0,1}
		=0.
	\end{align*}
	\noindent\textbf{Second identity ($k=1$):}
	Again only $\ell=1$ contributes, because $\tau^{2\ell-2}$ has degree at least two for $\ell\ge2$.
	Since $B_{0,1}=0$,
	\begin{align}\label{eq:first-defect-component}
		F_1(\tau)
		=\tau B_{1,1}+\tau^2B_{2,1}+\mathcal{O}(\tau^3).
	\end{align}
	Successive conjugation by the splitting factors gives
	\begin{align*}
		Z(\tau)B_{1,1}Z(-\tau)
		=B_{1,1}+\tau[A+B,B_{1,1}]+\mathcal{O}(\tau^2),
	\end{align*}
	so $Z(\tau)F_1(\tau)Z(-\tau)=\tau B_{1,1}+\mathcal{O}(\tau^2)$.
	Since $F_1(-\tau)=-\tau B_{1,1}+\mathcal{O}(\tau^2)$, the linear terms cancel:
	\begin{align*}
		K_1
		=[\tau]\big(Z(\tau)F_1(\tau)Z(-\tau)+F_1(-\tau)\big)
		=B_{1,1}-B_{1,1}
		=0.
	\end{align*}
	\noindent\textbf{Third identity ($k=2$):}
	For $n\ge2$, only $\ell=1,2$ contribute to the coefficient of degree two.
	For $\ell=1$, \eqref{eq:first-defect-component} and $Z(\tau)B_{2,1}Z(-\tau)=B_{2,1}+\mathcal{O}(\tau)$ give
	\begin{align*}
		Z(\tau)F_1(\tau)Z(-\tau)+F_1(-\tau)
		=\tau^2\bigl([A+B,B_{1,1}]+2B_{2,1}\bigr)+\mathcal{O}(\tau^3).
	\end{align*}
	For $\ell=2$, the prefactor $\tau^2$ means only $F_2(0)=B_{0,2}$ is needed:
	\begin{align*}
		\tau^2\bigl(Z(\tau)F_2(\tau)Z(-\tau)+F_2(-\tau)\bigr)
		=2\tau^2B_{0,2}+\mathcal{O}(\tau^3).
	\end{align*}
	Collecting these contributions finally yields
	\begin{align*}
		K_2
		&=[A+B,B_{1,1}]+2B_{2,1}+2B_{0,2}\\
		&=\left[A+B,\tfrac14[A,B]\right]-\frac18[A,[A,B]]-\frac18\big([A,[A,B]]-2[[A,B],B]\big)\\
		&=\frac14[A+B,[A,B]]-\frac14[A+B,[A,B]]
		=0,
	\end{align*}
	where the coefficients $B_{1,1},B_{2,1}$ and $B_{0,2}$ were taken from Remark~\ref{rem:defect-component-expansions}.
\end{example}
We conclude the section on the bounded case and now turn to the study of the local error in the setting of unbounded $A$ and $B$.
This step relies on the existence of $C_0$-groups generated by the Zassenhaus exponents $\mathcal{R}_{2\ell-1}(\tau)$, as well as on the truncation of the infinite series expansions of the local-error function.

\setcounter{equation}{0}
\section{Symmetric Zassenhaus splittings for unbounded operators}\label{sec:unbounded-zassenhaus}
\subsection{Preliminaries}
Let $A\colon\operatorname{Dom}(A)\to\mathcal{H}$ and $B\colon\operatorname{Dom}(B)\to\mathcal{H}$ be linear operators, where $\operatorname{Dom}(A)$ and $\operatorname{Dom}(B)$ are dense linear subspaces of $\mathcal{H}$.
Moreover, assume that a chosen realization $C$ of $A+B$, with domain $\operatorname{Dom}(C)\subseteq\mathcal{H}$, generates a strongly continuous semigroup on $\mathcal{H}$.
For $u_0\in\mathcal{H}$, consider the initial-value problem
\begin{align*}
	u'(t)
	=Cu(t),\qquad
	u(t_0)
	=u_0.
\end{align*}
Let $U(t)$, $t\ge0$, denote the semigroup generated by $C$.
The solution is
\begin{align*}
	u(t)
	=U(t-t_0)u_0,\qquad t\ge t_0.
\end{align*}
For $u_0\in\operatorname{Dom}(C)$, this solution is differentiable in $\mathcal{H}$ and satisfies $u'(t)=Cu(t)$; see \cite[Lemma~II.1.3]{Engel2000}.
As customary in the literature, the semigroup will be identified by
\begin{align*}
	\ee^{(A+B)(t-t_0)}
	:=U(t;t_0)=U(t-t_0).
\end{align*}
Here the exponential denotes the semigroup generated by $C$, rather than an operator power series; see e.g.\ \cite{Engel2000}.
\begin{example}[Linear Schrödinger equation]
	Let $\mathcal{H}=L^2(\BB{R}^d)$ and let $V\in L^\infty(\BB{R}^d;\BB{R})$.
	Set
	\begin{align*}
		A
		=\ii\Delta,\quad
		\operatorname{Dom}(A)
		=H^2(\BB{R}^d),\qquad
		B
		=-\ii V,\quad
		\operatorname{Dom}(B)
		=L^2(\BB{R}^d).	
	\end{align*}
	The operator $H=-\Delta+V$ is self-adjoint on the Sobolev space $H^2(\BB{R}^d)$; see e.g.\ \cite[Section~I.3.2]{Lubich2008}.
	Thus $C=-\ii H=A+B$, with $\operatorname{Dom}(C)=H^2(\BB{R}^d)$, generates the strongly continuous unitary group, with $U(t;t_0)=\ee^{-\ii(t-t_0)H}$.
	For $u_0\in L^2$, this gives the solution of
	\begin{align*}
		\ii u'(t)
		=(-\Delta+V)u(t),\qquad
		u(t_0)
		=u_0,
	\end{align*}
	and for $u_0\in H^2$, the solution is differentiable in $L^2$ and satisfies the equation there.
\end{example}
To extend the sZS to unbounded operators, we must give a meaning to each Zassenhaus factor and justify its differentiation.
We therefore assume that each $\mathcal{R}_{2\ell-1}(1)$ has a chosen realization generating a strongly continuous group, and interpret its exponential through that group.
For $1\le\ell\le n$, we interpret
\begin{align*}
	\ee^{\mathcal{R}_{2\ell-1}(\tau)}
	=\ee^{\tau^{2\ell-1}\mathcal{R}_{2\ell-1}(1)}	
\end{align*}
as this group evaluated at $t=\tau^{2\ell-1}$.
On the domain of $\mathcal{R}_{2\ell-1}(1)$, semigroup differentiation gives
\begin{align*}
	\frac{\mathrm{d}}{\mathrm{d}\tau}\ee^{\mathcal{R}_{2\ell-1}(\tau)}\psi
	=(2\ell-1)\tau^{2\ell-2}\mathcal{R}_{2\ell-1}(1)\ee^{\mathcal{R}_{2\ell-1}(\tau)}\psi.
\end{align*}
Recall that for $\ell=0$, the corresponding factor is $\ee^{\mathcal{R}_{-1}(\tau)}=\ee^{\tau A/2}$.

A naive approach to proving the local-error order for unbounded operators is to apply Taylor's theorem to both $U(\tau)u_0$ and $Z_{2n-1}(\tau)u_0$ and verify cancellation up to degree $2n$.
However, the Taylor remainder of the exact solution takes the form
\begin{align*}
	\frac{1}{(2n)!}\int_0^\tau(\tau-s)^{2n}U(s)C^{2n+1}u_0\,\mathrm{d}s.
\end{align*}
This argument therefore requires $u_0\in\operatorname{Dom}(C^{2n+1})$, which can impose substantially stronger regularity than estimates expressed in terms of commutators.
Moreover, the approach pursued in Theorem~\ref{thm:bounded-local-error} is not applicable, since it relies on the sBCH formula presented in \eqref{eq:sbch-series}, which is not necessarily convergent in the unbounded setting.\footnote{The truncated sBCH expression combines unbounded commutators of different lengths, multiplied by different powers of $\tau$.
	Defining its exponential requires additional domain and generation assumptions, which we avoid in the following analysis.}
To establish a local-error bound for unbounded operators, we instead exploit the integral representation of the error given in \eqref{eq:duhamel-error-representation}, together with the expansion of $G(\tau)$ in \eqref{eq:factored-defect-series}.

In the unbounded setting, remainder estimates must specify the spaces and norms involved.
We therefore introduce the following notation:
\begin{definition}
	For normed spaces $\mathcal{X}$ and $\mathcal{Y}$, we denote by $\mathcal{L}(\mathcal{X},\mathcal{Y})$ the space of bounded linear operators from $\mathcal{X}$ into $\mathcal{Y}$.
	Let $\mathcal{H}$ be a Hilbert space and let $\mathcal{D}$ be a Banach space densely and continuously embedded into $\mathcal{H}$.
	For a family $C(\tau)\in\mathcal{L}(\mathcal{D},\mathcal{H})$, $\tau\ge 0$, and $k\ge 0$, we write $C(\tau)=\mathcal{O}(\tau^k)$, if there exist constants $c>0$ and $\tau_0>0$ such that
	\begin{align*}
		\|C(\tau)\psi\|_{\mathcal{H}}
		\le c\tau^k\|\psi\|_{\mathcal{D}}\quad
		\text{for all $\psi\in\mathcal{D}$ and all $0\leq\tau\leq\tau_0$.}
	\end{align*}
	Moreover, throughout the unbounded-operator analysis, operator-valued integrals are understood after application to each $\psi\in\mathcal{D}$ as Bochner integrals in $\mathcal{H}$.
\end{definition}
Obtaining a bound for the factored defect requires a generalization of Proposition~\ref{prop:bounded-conjugation-expansion}, which we present in the following subsection.
This generalization, in which the remainder terms are controlled explicitly, will in turn lead to regularity restrictions of the initial condition.

\subsection{Finite conjugation expansions and overshoots}
Before presenting the analogue of Proposition~\ref{prop:bounded-conjugation-expansion} for unbounded operators, we generalize the expansion of
\begin{align}\label{eq:single-conjugation-definition}
	H(t)
	:=\ee^{t^kX}Y\ee^{-t^kX},
\end{align}
where $X$ and $Y$ are homogeneous Lie polynomials in $A,B$, and $k=\ell(X)\ge1$.
We work under Assumption~\ref{ass:unbounded-framework}, stated below, with commutator depth at least $\nu+1$ for the Taylor expansion below; identities and derivatives are understood after application to $\psi\in\mathcal{D}$.
Applying Taylor's theorem\footnote{For unbounded operators, this expansion can be rigorously justified by iterating the fundamental theorem of calculus.} gives
\begin{align}\label{eq:conjugation-taylor-expansion}
	H(t)
	=\sum_{s=0}^{\nu}\frac{t^s}{s!}H^{(s)}(0)+\int_0^t\frac{(t-\tau)^\nu}{\nu!}H^{(\nu+1)}(\tau)\,\mathrm{d}\tau.
\end{align}
Unlike in the bounded case, the integral term is not yet in the desired form.
To see this, we first derive an explicit formula for the derivatives of $H$:
\begin{lemma}\label{lem:conjugation-derivatives}
	For $1\le s\le J$, the $s$-th derivative of the function $H$ in \eqref{eq:single-conjugation-definition} is given by
	\begin{align}\label{eq:conjugation-derivatives}
		H^{(s)}(t)
		=\sum_{m\in\mathcal{M}_{k,s}}c_{m,s}\ee^{t^kX}\operatorname{ad}_X^{|m|}(Y)\ee^{-t^kX}t^{k|m|-s},
	\end{align}
	with positive constants $c_{m,s}>0$ given in \eqref{eq:derivative-combinatorial-coefficients}, and
	\begin{align}\label{eq:derivative-multiindices}
		\mathcal{M}_{k,s}
		:=\{m\in\BB{Z}_{\ge 0}^s\,\colon\,\sum_{j=1}^sjm_j=s,\,\text{$m_j=0$ for all $j>k$}\}.
	\end{align}
	In particular, $H^{(s)}(0)$ vanishes unless $s\equiv 0\pmod k$.
	In that case, set $j=s/k$.
	Then
	\begin{align}\label{eq:conjugation-derivatives-at-zero}
		H^{(s)}(0)
		=\frac{s!}{j!}\operatorname{ad}_X^j(Y).
	\end{align}
\end{lemma}
\begin{proof}
	See Appendix~\ref{app:conjugation-derivatives}.
\end{proof}

According to the previous lemma, the derivative $H^{(\nu+1)}$ appearing in \eqref{eq:conjugation-taylor-expansion} naturally decomposes into a finite sum, with one contribution for each $m\in\mathcal{M}_{k,\nu+1}$.
Moreover, inserting the expression for the derivative produces the summands
\begin{align}\label{eq:taylor-remainder-contributions}
	\int_0^t\frac{(t-\tau)^\nu}{\nu!}c_{m,\nu+1}\ee^{\tau^kX}\operatorname{ad}_X^{|m|}(Y)\ee^{-\tau^kX}\tau^{k|m|-(\nu+1)}\,\mathrm{d}\tau
	=\mathcal{O}\big(t^{k|m|}\big),
\end{align}
and the corresponding indices of the sum satisfy
\begin{align}\label{eq:multiindex-length-bound}
	m_*
	:=\left\lfloor\frac{\nu}{k}\right\rfloor+1
	\le |m|
	\le\nu+1
	=:m^*.
\end{align}
Hence the individual contributions range from\footnote{This follows from $k(\lfloor\nu/k\rfloor+1)\ge\nu+1$.} order $\mathcal{O}(t^{km_*})=\mathcal{O}(t^{\nu+1})$ to order $\mathcal{O}(t^{k(\nu+1)})$.
The latter contribution is precisely the one involving the highest-order commutator $\operatorname{ad}_X^{\nu+1}(Y)$.

Since estimates involving $\operatorname{ad}_X^{\nu+1}(Y)$ typically require substantially stronger regularity assumptions than those involving lower-order commutators, it is advantageous to isolate this contribution from all terms with $|m|\le\nu$.
The next lemma shows that the lower-order contributions can be extracted explicitly, leaving a remainder controlled by the highest-order commutator and the local group bounds.

To ensure that the iterated commutators appearing below are well defined and that the corresponding remainder estimates hold uniformly on a common regularity space, we adopt the following functional-analytic framework throughout the unbounded-operator analysis.
\begin{assumption}[Regularity framework]\label{ass:unbounded-framework}
	Let $X\colon\operatorname{Dom}(X)\subseteq\mathcal{H}\to\mathcal{H}$ and $Y\colon\operatorname{Dom}(Y)\subseteq\mathcal{H}\to\mathcal{H}$ be linear operators.
	Moreover, let $\mathcal{D}$ be a Banach space densely and continuously embedded into $\mathcal{H}$, and let $\mathcal{D}_0$ be a dense linear subspace of $\mathcal{D}$.
	For a required commutator depth $J\ge1$, we assume:
	\begin{enumerate}
		\item\textbf{Controlled commutators.}
			All products defining $\operatorname{ad}_X^j(Y)$, $0\le j\le J$, are defined on $\mathcal{D}_0$, and these commutators extend boundedly from $\mathcal{D}$ to $\mathcal{H}$.
		\item\textbf{Regularity-preserving flow.}
			The operator $X$ generates a strongly continuous \emph{group} on $\mathcal{H}$ which preserves $\mathcal{D}$ and restricts to a strongly continuous \emph{group} on $\mathcal{D}$.
		\item\textbf{Integral conjugation identity.}
			For $0\le j<J$, $t\in\BB{R}$, and $\psi\in\mathcal{D}$,
			\begin{align*}
				\ee^{tX}\operatorname{ad}_X^j(Y)\ee^{-tX}\psi
				=\operatorname{ad}_X^j(Y)\psi+\int_0^t\ee^{sX}\operatorname{ad}_X^{j+1}(Y)\ee^{-sX}\psi\,\mathrm{d}s.
			\end{align*}
	\end{enumerate}
\end{assumption}
We illustrate these assumptions with a concrete example from the Schr\"odinger equation.
\begin{example}[Schr\"odinger operators]
	Let $\mathcal{H}=L^2(\BB{R}^d)$, $X=\ii\Delta$, and $Y=-\ii V$, where $V\in C^\infty(\BB{R}^d;\BB{R})$.
	We choose $V$ and $\mathcal{D}$ so that the three conditions of Assumption~\ref{ass:unbounded-framework} hold for this pair and the required $J\ge 1$.
	Note that boundedness of $V$ is not required; unbounded potentials may require a weighted regularity space.
	For example, the quadratic potential $V(x)=|x|^2$ satisfies these conditions for every $J\ge1$ with
	\begin{align*}
		\mathcal{D}
		=\{\psi\in H^2(\BB{R}^d)\,\colon\,|x|^2\psi\in L^2(\BB{R}^d)\},\qquad
		\|\psi\|_{\mathcal{D}}^2
		=\|\psi\|_{H^2}^2+\||x|^2\psi\|_{L^2}^2.
	\end{align*}
	This is a weighted Sobolev space.\footnote{It is denoted by $\Sigma^2$ in \cite[Remark~7.1, with $\varepsilon=1$]{Lasser2020} and by $M_2^2$, the \emph{Shubin--Sobolev space} of order two, in \cite[Section~2]{Wahlberg2022}.
	The corresponding norms are equivalent to $\|\cdot\|_{\mathcal{D}}$.}
	We take $\mathcal{D}_0=C_c^\infty(\BB R^d)$.
	\begin{enumerate}
		\item\textbf{Controlled commutators.}
			The space $\mathcal{D}$ is densely and continuously embedded into $L^2$, and $\mathcal{D}_0$ is dense in $\mathcal{D}$.
			All required products are defined on $\mathcal{D}_0$, and
			\begin{align*}
				[X,Y]
				=4x\cdot\nabla+2d\operatorname{Id},\qquad
				\operatorname{ad}_X^2(Y)
				=8\ii\Delta,\qquad
				\operatorname{ad}_X^j(Y)
				=0\quad\text{for $j\ge 3$}.
			\end{align*}
			These operators and $Y$ are bounded from $\mathcal{D}$ to $L^2$; the mixed terms are controlled by
			\begin{align*}
				\|x_r\partial_{x_q}\psi\|_{L^2}
				\le c\|\psi\|_{\mathcal{D}},\quad
				r,q\in\{1,\ldots,d\},
			\end{align*}
			which can be obtained by integration by parts.
		\item\textbf{Regularity-preserving flow.}
			The operator $X=\ii\Delta$, with domain $H^2(\BB{R}^d)$, generates the strongly continuous unitary group $\ee^{\ii t\Delta}$ on $\mathcal{H}=L^2(\BB{R}^d)$.
			This group preserves $\mathcal{D}$ and restricts to a strongly continuous group on $\mathcal{D}$.
			Indeed, it commutes with derivatives, and
			\begin{align}\label{eq:free-flow-coordinate-identity}
				x_r\ee^{\ii t\Delta}\psi
				=\ee^{\ii t\Delta}(x_r-2\ii t\partial_{x_r})\psi,\quad
				1\le r\le d,
			\end{align}
			applied twice, controls the quadratic weight $|x|^2$ using the preceding mixed-term estimate.
		\item\textbf{Integral conjugation identity.}
			For $t\in\BB{R}$ and $\psi\in\mathcal{D}$, the coordinate identity \eqref{eq:free-flow-coordinate-identity} gives
			\begin{align*}
				\ee^{tX}Y\ee^{-tX}\psi
				&=Y\psi+t[X,Y]\psi+\frac{t^2}{2}\operatorname{ad}_X^2(Y)\psi,\\
				\ee^{tX}[X,Y]\ee^{-tX}\psi
				&=[X,Y]\psi+t\operatorname{ad}_X^2(Y)\psi.
			\end{align*}
			Moreover, $\operatorname{ad}_X^2(Y)=8\ii\Delta$ commutes with the free group, and all higher commutators vanish.
			Integrating the conjugation formula for each next commutator therefore yields, for every $0\le j<J$,
			\begin{align*}
				\ee^{tX}\operatorname{ad}_X^j(Y)\ee^{-tX}\psi
				=\operatorname{ad}_X^j(Y)\psi+\int_0^t\ee^{sX}\operatorname{ad}_X^{j+1}(Y)\ee^{-sX}\psi\,\mathrm{d}s.
			\end{align*}
	\end{enumerate}
\end{example}
Under the above assumptions, we can now extract the lower-order terms from the integral remainder in \eqref{eq:taylor-remainder-contributions}, leaving a remainder controlled by the highest-order commutator and the local group bounds.
\begin{lemma}\label{lem:integral-remainder-decomposition}
	Let $X$ and $Y$ satisfy Assumption~\ref{ass:unbounded-framework} with maximal commutator depth $J=\nu+1$, where $\nu\ge 0$, and let $k=\ell(X)$.
	Then, for all $t\ge0$, the integral remainder in \eqref{eq:conjugation-taylor-expansion} admits the decomposition
	\begin{align*}
		\int_0^t\frac{(t-\tau)^\nu}{\nu!}H^{(\nu+1)}(\tau)\,\mathrm{d}\tau
		=\sum_{j=\lfloor\nu/k\rfloor+1}^{\nu}\frac{t^{jk}}{j!}\operatorname{ad}_X^j(Y)+R_{\nu+1}(t),
	\end{align*}
	as an identity in $\mathcal{L}(\mathcal{D},\mathcal{H})$, where $R_{\nu+1}(t)=\mathcal{O}(t^{k(\nu+1)})$.
	In particular, the remainder bound depends only on $k$, $\nu$, the $\mathcal{L}(\mathcal{D},\mathcal{H})$-norm of $\operatorname{ad}_X^{\nu+1}(Y)$, and the local group bounds in Assumption~\ref{ass:unbounded-framework}.
\end{lemma}
\begin{proof}
	For convenience, introduce the family of conjugated operators
	\begin{align*}
		H_j(\tau)
		:=\ee^{\tau^kX}\operatorname{ad}_X^j(Y)\ee^{-\tau^kX},\quad
		0
		\le j
		\le\nu+1.
	\end{align*}
	By Assumption~\ref{ass:unbounded-framework}, every operator $\operatorname{ad}_X^j(Y)$ admits a bounded extension from $\mathcal{D}$ into $\mathcal{H}$.
	Moreover, the group generated by $X$ is locally uniformly bounded on both $\mathcal{H}$ and $\mathcal{D}$.
	Consequently, $H_j(\tau)\in\mathcal{L}(\mathcal{D},\mathcal{H})$ for every $\tau$ in a bounded interval.
	The integral conjugation identity in Assumption~\ref{ass:unbounded-framework} gives
	\begin{align}\label{eq:commutator-conjugation-integral}
		H_j(\tau)
		=\operatorname{ad}_X^j(Y)+\int_0^\tau ks^{k-1}H_{j+1}(s)\,\mathrm{d}s,\quad
		0
		\le j
		\le\nu,
	\end{align}
	as an identity in $\mathcal{L}(\mathcal{D},\mathcal{H})$.
	For $j=\nu$, formula \eqref{eq:commutator-conjugation-integral} reads
	\begin{align*}
		H_\nu(\tau)
		=\operatorname{ad}_X^\nu(Y)+\int_0^\tau ks^{k-1}H_{\nu+1}(s)\,\mathrm{d}s.
	\end{align*}
	Iterating \eqref{eq:commutator-conjugation-integral} yields, for $1\le j\le\nu+1$,
	\begin{align}\label{eq:commutator-conjugation-decomposition}
		H_j(\tau)
		=P_j(\tau)+E_j(\tau),
	\end{align}
	where we introduced the two contributions
	\begin{align*}
		P_j(\tau)
		&=\sum_{\ell=0}^{\nu-j}\frac{\tau^{\ell k}}{\ell!}\operatorname{ad}_X^{j+\ell}(Y),\,1\le j\le\nu,\quad
		P_{\nu+1}
		=0,\\
		E_j(\tau)
		&=\int_{\Delta^{\nu+1-j}(\tau)}\hspace{-5mm}H_{\nu+1}(s_{\nu+1-j})\prod_{i=1}^{\nu+1-j}ks_i^{k-1}\,\mathrm{d}V_{\nu+1-j},\,1\le j\le\nu,\quad
		E_{\nu+1}
		=H_{\nu+1},
	\end{align*}
	and use the shorthand notation
	\begin{align*}
		\int_{\Delta^L(\tau)}f(s_1,\ldots,s_L)\,\mathrm{d}V_L
		:=\int_0^\tau\int_0^{s_1}\cdots\int_0^{s_{L-1}}f(s_1,\ldots,s_L)\,\mathrm{d}s_L\cdots\mathrm{d}s_1.
	\end{align*}
	In particular, all identities in \eqref{eq:commutator-conjugation-decomposition} hold in $\mathcal{L}(\mathcal{D},\mathcal{H})$.
	For $j=\nu+1$, subsequent simplex expressions are understood without an inner integral, with $s_0=\tau$ and the empty product equal to one.
	
	Substituting the representation of $H^{(\nu+1)}$ from \eqref{eq:conjugation-derivatives} into the integral remainder in \eqref{eq:conjugation-taylor-expansion}, we obtain
	\begin{align*}
		\int_0^t\frac{(t-\tau)^\nu}{\nu!}H^{(\nu+1)}(\tau)\,\mathrm{d}\tau
		&=\int_0^t\frac{(t-\tau)^\nu}{\nu!}\sum_{m\in\mathcal{M}_{k,\nu+1}}c_{m,\nu+1}H_{|m|}(\tau)\tau^{k|m|-(\nu+1)}\,\mathrm{d}\tau\\
		&=\sum_{j=1}^{\nu+1}\sum_{\begin{subarray}{c}m\in\mathcal{M}_{k,\nu+1}\\|m|=j\end{subarray}}c_{m,\nu+1}\int_0^t\frac{(t-\tau)^\nu}{\nu!}H_j(\tau)\tau^{kj-(\nu+1)}\,\mathrm{d}\tau.
	\end{align*}
	Inserting the decomposition \eqref{eq:commutator-conjugation-decomposition} then gives
	\begin{align*}
		\int_0^t\frac{(t-\tau)^\nu}{\nu!}H^{(\nu+1)}(\tau)\,\mathrm{d}\tau
		=P(t)+R_{\nu+1}(t),
	\end{align*}
	where
	\begin{align*}
		P(t)
		&=\sum_{j=1}^{\nu}\sum_{\begin{subarray}{c}m\in\mathcal{M}_{k,\nu+1}\\|m|=j\end{subarray}}c_{m,\nu+1}\int_0^t\frac{(t-\tau)^\nu}{\nu!}P_j(\tau)\tau^{kj-(\nu+1)}\,\mathrm{d}\tau,\\
		R_{\nu+1}(t)
		&=\sum_{j=1}^{\nu+1}\sum_{\begin{subarray}{c}m\in\mathcal{M}_{k,\nu+1}\\|m|=j\end{subarray}}c_{m,\nu+1}\int_0^t\frac{(t-\tau)^\nu}{\nu!}E_j(\tau)\tau^{kj-(\nu+1)}\,\mathrm{d}\tau.
	\end{align*}
	We first determine the polynomial contribution $P(t)$.
	According to the bounds for $|m|$ presented in \eqref{eq:multiindex-length-bound}, we conclude that $\left\{m\in\mathcal{M}_{k,\nu+1}\,\colon\,|m|=j\right\}=\emptyset$ for $j\le\lfloor\nu/k\rfloor$.
	Moreover,
	\begin{align*}
		\int_0^t\frac{(t-\tau)^\nu}{\nu!}\left(\sum_{\ell=0}^{\nu-j}\operatorname{ad}_X^{j+\ell}(Y)\frac{\tau^{\ell k}}{\ell!}\right)\tau^{kj-(\nu+1)}\,\mathrm{d}\tau
		=\sum_{\ell=0}^{\nu-j}\operatorname{ad}_X^{j+\ell}(Y)\frac{t^{k(j+\ell)}}{\ell!\big[k(j+\ell)\big]_{\nu+1}},
	\end{align*}
	where $[k]_j=k(k-1)\cdots(k-j+1)$ denotes the descending factorial.
	It follows that
	\begin{align*}
		P(t)
		&=\sum_{j=\lfloor\nu/k\rfloor+1}^{\nu}\sum_{\ell=0}^{\nu-j}\operatorname{ad}_X^{j+\ell}(Y)\frac{t^{k(j+\ell)}}{\ell!\bigl[k(j+\ell)\bigr]_{\nu+1}}\sum_{\begin{subarray}{c}m\in\mathcal{M}_{k,\nu+1}\\|m|=j\end{subarray}}c_{m,\nu+1}\\
		&=\sum_{j=\lfloor\nu/k\rfloor+1}^{\nu}s_{j,\nu+1}\sum_{r=j}^{\nu}\operatorname{ad}_X^r(Y)\frac{t^{kr}}{(r-j)![kr]_{\nu+1}}\\
		&=\sum_{r=\lfloor\nu/k\rfloor+1}^{\nu}\operatorname{ad}_X^r(Y)\frac{t^{kr}}{[kr]_{\nu+1}}\sum_{j=1}^{r}\frac{s_{j,\nu+1}}{(r-j)!},
	\end{align*}
	with the constants
	\begin{align*}
		s_{j,\nu+1}
		:=\sum_{\begin{subarray}{c}m\in\mathcal{M}_{k,\nu+1}\\|m|=j\end{subarray}}c_{m,\nu+1}.
	\end{align*}
	Using the combinatorial identity
	\begin{align}\label{eq:combinatorial-identity}
		\sum_{j=1}^{r}\frac{s_{j,\nu+1}}{(r-j)!}
		=\frac{[kr]_{\nu+1}}{r!},
	\end{align}
which is proven in Appendix~\ref{app:combinatorial-identity}, we conclude that
	\begin{align}\label{eq:polynomial-contribution}
		P(t)
		=\sum_{r=\lfloor\nu/k\rfloor+1}^{\nu}\frac{t^{kr}}{r!}\operatorname{ad}_X^r(Y).
	\end{align}
	It therefore remains to estimate the contribution from $R_{\nu+1}(t)$.
	
	Fix $T>0$ and let $0\le t\le T$.
	Strong continuity of the group on $\mathcal{H}$ and $\mathcal{D}$ implies local uniform boundedness; see \cite[Proposition~I.5.5]{Engel2000}, applied in both time directions.
	Thus there exist constants $m_{\mathcal{H},T},m_{\mathcal{D},T}>0$ such that
	\begin{align*}
		\sup_{|s|\leq T^k}\|\ee^{sX}\|_{\mathcal{L}(\mathcal{H})}
		\le m_{\mathcal{H},T},\qquad
		\sup_{|s|\leq T^k}\|\ee^{sX}\|_{\mathcal{L}(\mathcal{D})}
		\le m_{\mathcal{D},T}.
	\end{align*}
	Hence, for $0\le s\le t$ and $\psi\in\mathcal{D}$,
	\begin{align*}
		\|H_{\nu+1}(s)\psi\|_{\mathcal{H}}
		\le m_{\mathcal{H},T}\left\|\operatorname{ad}_X^{\nu+1}(Y)\right\|_{\mathcal{L}(\mathcal{D},\mathcal{H})}m_{\mathcal{D},T}\|\psi\|_{\mathcal{D}}.
	\end{align*}
	Thus, there exists a constant $m_T>0$, independent of $(s,t,\psi)$, such that
	\begin{align}\label{eq:highest-commutator-bound}
		\sup_{0\leq s\leq t}\|H_{\nu+1}(s)\psi\|_{\mathcal{H}}
		\le m_T\|\psi\|_{\mathcal{D}}.
	\end{align}
	Using the triangle inequality and the definition of $E_j$, we then obtain
	\begin{align*}
		\|R_{\nu+1}(t)\psi\|_{\mathcal{H}}
		&\le\sum_{j=1}^{\nu+1}\sum_{\begin{subarray}{c}m\in\mathcal{M}_{k,\nu+1}\\|m|=j\end{subarray}}c_{m,\nu+1}\int_0^t\frac{(t-\tau)^\nu}{\nu!}\\
		&\qquad\times\left(\int_{\Delta^{\nu+1-j}(\tau)}\|H_{\nu+1}(s_{\nu+1-j})\psi\|_{\mathcal{H}}\prod_{i=1}^{\nu+1-j}ks_i^{k-1}\,\mathrm{d}V_{\nu+1-j}\right)\tau^{kj-(\nu+1)}\,\mathrm{d}\tau.
	\end{align*}
	By the bound in \eqref{eq:highest-commutator-bound} and the simple formula
	\begin{align*}
		\int_{\Delta^{\nu+1-j}(\tau)}\prod_{i=1}^{\nu+1-j}ks_i^{k-1}\,\mathrm{d}V_{\nu+1-j}
		=\frac{\tau^{k(\nu+1-j)}}{(\nu+1-j)!},
	\end{align*}
	it therefore follows that
	\begin{align*}
		\|R_{\nu+1}(t)\psi\|_{\mathcal{H}}
		\le m_T\|\psi\|_{\mathcal{D}}\sum_{j=1}^{\nu+1}\sum_{\begin{subarray}{c}m\in\mathcal{M}_{k,\nu+1}\\|m|=j\end{subarray}}c_{m,\nu+1}\int_0^t\frac{(t-\tau)^\nu}{\nu!(\nu+1-j)!}\tau^{(k-1)(\nu+1)}\,\mathrm{d}\tau.
	\end{align*}
	The remaining integral satisfies
	\begin{align*}
		\int_0^t\frac{(t-\tau)^\nu}{\nu!(\nu+1-j)!}\tau^{(k-1)(\nu+1)}\,\mathrm{d}\tau
		=\frac{t^{k(\nu+1)}}{(\nu+1-j)!}\frac{1}{\bigl[k(\nu+1)\bigr]_{\nu+1}}.
	\end{align*}
	Consequently,
	\begin{align*}
		\|R_{\nu+1}(t)\psi\|_{\mathcal{H}}
		\le m_T\|\psi\|_{\mathcal{D}}\frac{t^{k(\nu+1)}}{\bigl[k(\nu+1)\bigr]_{\nu+1}}\sum_{j=1}^{\nu+1}\frac{s_{j,\nu+1}}{(\nu+1-j)!}.
	\end{align*}
	Applying once more the combinatorial identity \eqref{eq:combinatorial-identity} with $r=\nu+1$ finally gives
	\begin{align*}
		\|R_{\nu+1}(t)\|_{\mathcal{L}(\mathcal{D},\mathcal{H})}
		\le\frac{m_T}{(\nu+1)!}t^{k(\nu+1)}.
	\end{align*}
	Combining this estimate with \eqref{eq:polynomial-contribution} proves the result.
\end{proof}

Lemma~\ref{lem:conjugation-derivatives} identifies the Taylor coefficients of $\ee^{t^kX}Y\ee^{-t^kX}$, while Lemma~\ref{lem:integral-remainder-decomposition} extracts the lower-order terms from the integral remainder.
Together, they yield a finite commutator expansion whose remainder is controlled solely by the next commutator:
\begin{proposition}\label{prop:single-conjugation-expansion}
	Under the assumptions of Lemma~\ref{lem:integral-remainder-decomposition}, we have, for all $t\ge 0$,
	\begin{align*}
		\ee^{t^kX}Y\ee^{-t^kX}
		=\sum_{j=0}^{\nu}\frac{t^{jk}}{j!}\operatorname{ad}_X^j(Y)+R_{\nu+1}(t)
	\end{align*}
	as an identity in $\mathcal{L}(\mathcal{D},\mathcal{H})$, where $R_{\nu+1}(t)=\mathcal{O}(t^{k(\nu+1)})$.
	In particular, the remainder bound depends only on $k$, $\nu$, the $\mathcal{L}(\mathcal{D},\mathcal{H})$-norm of $\operatorname{ad}_X^{\nu+1}(Y)$, and the local group bounds in Assumption~\ref{ass:unbounded-framework}.
\end{proposition}
\begin{proof}
	For $\nu=0$, the integral conjugation identity gives directly
	\begin{align*}
		\ee^{t^kX}Y\ee^{-t^kX}-Y
		=\int_0^{t^k}\ee^{sX}[X,Y]\ee^{-sX}\,\mathrm{d}s
		=\mathcal{O}_{\mathcal{L}(\mathcal{D},\mathcal{H})}(t^k).
	\end{align*}
	Hence it remains to consider $\nu\ge1$.
	
	Let $\psi\in\mathcal{D}$.
	By the formula \eqref{eq:conjugation-derivatives-at-zero}, the derivatives of $H$ at zero satisfy
	\begin{align*}
		H^{(s)}(0)\psi
		=\mathbf{1}_{k\scalebox{0.65}{$\BB{Z}$}}(s)\frac{s!}{(s/k)!}\operatorname{ad}_X^{s/k}(Y)\psi,
	\end{align*}
	where $\mathbf{1}_{k\scalebox{0.65}{\BB{Z}}}\colon\BB{Z}\to\{0,1\}$ denotes the indicator function on integer multiples of $k$.
	Note that this expression is understood to be zero when $s$ is not a multiple of $k$, without evaluating the factorial.
	Consequently, applying Taylor's formula to the $\mathcal{H}$-valued function $H(t)\psi$, we obtain
	\begin{align*}
		\ee^{t^kX}Y\ee^{-t^kX}\psi
		&=\sum_{s=0}^{\nu}\frac{t^s}{s!}H^{(s)}(0)\psi+\int_0^t\frac{(t-\tau)^\nu}{\nu!}H^{(\nu+1)}(\tau)\psi\,\mathrm{d}\tau\\
		&=\sum_{s=0}^{\nu}\frac{t^s}{s!}\mathbf{1}_{k\scalebox{0.65}{$\BB{Z}$}}(s)\frac{s!}{(s/k)!}\operatorname{ad}_X^{s/k}(Y)\psi+\int_0^t\frac{(t-\tau)^\nu}{\nu!}H^{(\nu+1)}(\tau)\psi\,\mathrm{d}\tau\\
		&=\sum_{j=0}^{\lfloor\nu/k\rfloor}\frac{t^{jk}}{j!}\operatorname{ad}_X^j(Y)\psi+\int_0^t\frac{(t-\tau)^\nu}{\nu!}H^{(\nu+1)}(\tau)\psi\,\mathrm{d}\tau.
	\end{align*}
	By Lemma~\ref{lem:integral-remainder-decomposition}, the integral remainder admits the decomposition
	\begin{align*}
		\int_0^t\frac{(t-\tau)^\nu}{\nu!}H^{(\nu+1)}(\tau)\,\mathrm{d}\tau
		=\sum_{j=\lfloor\nu/k\rfloor+1}^{\nu}\frac{t^{jk}}{j!}\operatorname{ad}_X^j(Y)+R_{\nu+1}(t),
	\end{align*}
	where $R_{\nu+1}(t)=\mathcal{O}(t^{k(\nu+1)})$.
	Combining the two finite sums therefore gives
	\begin{align*}
		\ee^{t^kX}Y\ee^{-t^kX}\psi
		=\sum_{j=0}^{\nu}\frac{t^{jk}}{j!}\operatorname{ad}_X^j(Y)\psi+R_{\nu+1}(t)\psi.
	\end{align*}
	Since this identity holds for every $\psi\in\mathcal{D}$, and all operators belong to $\mathcal{L}(\mathcal{D},\mathcal{H})$, it holds in $\mathcal{L}(\mathcal{D},\mathcal{H})$.
\end{proof}

We now extend Proposition~\ref{prop:bounded-conjugation-expansion} to unbounded operators.
Unlike in the bounded setting, the existence of the expansion requires additional regularity assumptions.
Let us first introduce the following two definitions:
\begin{definition}
	Let $\GG{L}(A,B)$ denote the free Lie algebra generated by the letters $A$ and $B$.
	For an integer $q\ge1$, let $\GG{L}_q(A,B)$ be the span of nested commutators of length $q$, where length counts the occurrences of $A$ and $B$.
	For a nonzero homogeneous Lie polynomial $F\in\GG{L}_q(A,B)$, write $\ell(F)=q$.
	For an integer $P^*\ge1$, define
	\begin{align*}
		\GG{L}_{\leq P^*}(A,B)
		:=\operatorname{span}\{F\,\colon\,\text{$F$ is a nested commutator in $A$, $B$ of length at most $P^*$}\}.
	\end{align*}
\end{definition}
\begin{definition}
	For a given integer $P\ge 1$, we denote by
	\begin{align*}
		\langle P\rangle_k
		:=P\bmod k
		=P-k\left\lfloor\frac Pk\right\rfloor
	\end{align*}
	the distance\,\footnote{Note that $0\le\langle P\rangle_k\le k-1$.} from $P$ to the largest multiple of $k\ge1$ not exceeding $P$.
\end{definition}
We now state the finite expansion for successive conjugations by unbounded operators:
\begin{theorem}\label{thm:multiple-conjugation-expansion}
	Let $n\ge1$, and let $X_1,\ldots,X_n$ and $Y$ be nonzero homogeneous Lie polynomials in $A,B$, with $k_j:=\ell(X_j)$, $1\le j\le n$.
	For an integer $P\ge1$, set
	\begin{align*}
		k_*
		:=\max_{1\le j\le n}k_j,\qquad
		P^*
		:=P+\ell(Y)+k_*.	
	\end{align*}
	Let $\mathcal{D}$ and $\mathcal{D}_0$ be as in Assumption~\ref{ass:unbounded-framework}.
	For the successive conjugations, we require the following versions of its three conditions, all on the same spaces $\mathcal{D}_0\subseteq\mathcal{D}\subseteq\mathcal{H}$:
	\begin{enumerate}
		\item All products defining $F\in\GG L_{\le P^*}(A,B)$ are defined on $\mathcal{D}_0$, and every such $F$ admits a bounded extension from $\mathcal{D}$ to $\mathcal{H}$.
		\item A chosen realization of each $X_j$ generates a strongly continuous group on $\mathcal{H}$ preserving $\mathcal{D}$ and restricting to a strongly continuous group on $\mathcal{D}$.
			These realizations agree with their algebraic expressions on $\mathcal{D}_0$.
		\item For $1\le j\le n$, the integral conjugation identity in Assumption~\ref{ass:unbounded-framework} holds with $X=X_j$ and $Y=F$ for every $F\in\GG{L}_q(A,B)$ with $1\le q\le P^*-k_j$.
	\end{enumerate}
	Then there exist coefficients $B_{p,n}\in\GG{L}_{p+\ell(Y)}(A,B)$, $0\le p\le P$, such that
	\begin{align*}
		A_n(t)
		=\left(\ee^{t^{k_n}X_n}\cdots\ee^{t^{k_1}X_1}\right)Y\left(\ee^{-t^{k_1}X_1}\cdots\ee^{-t^{k_n}X_n}\right)
		=\sum_{p=0}^{P}B_{p,n}t^p+R_{P+1,n}(t)
	\end{align*}
	as an identity in $\mathcal{L}(\mathcal{D},\mathcal{H})$, with
	\begin{align*}
		\|R_{P+1,n}(t)\|_{\mathcal{L}(\mathcal{D},\mathcal{H})}
		\le c_Tt^{P+1},\quad
		0\le t\le T.	
	\end{align*}
	For each fixed $T>0$, the constant $c_T$ depends only on $P$, the lengths $k_j,\ell(Y)$, the local group bounds, and the $\mathcal{L}(\mathcal{D},\mathcal{H})$-norms of the Lie polynomials occurring in the expansion and remainder.
\end{theorem}
\begin{remark}
	The key distinction between Proposition~\ref{prop:bounded-conjugation-expansion} (bounded setting) and Theorem~\ref{thm:multiple-conjugation-expansion} (unbounded setting) lies in the type of expansion that can be justified.
	When $A$ and $B$ are bounded, the full infinite expansion \eqref{eq:bounded-conjugation-expansion} is valid, and it is natural -- indeed harmless -- to allow compositions and commutators of arbitrary length: every composition remains a bounded operator, so no domain restrictions on elements of the space ever arise and convergence of the series is not an issue.
	
	In the unbounded setting of Theorem~\ref{thm:multiple-conjugation-expansion}, however, the infinite expansion is not justified under the stated assumptions: we can only assert a finite expansion up to a given order $P\ge 1$, with a remainder $\mathcal{O}(t^{P+1})$.
	The finite-expansion argument gives the sufficient length bound $P^*=P+\ell(Y)+k_*$.
	When $k_*>1$, this bound exceeds $P+\ell(Y)+1$, the length suggested by the first omitted degree, and we call this excess an `overshoot'.
\end{remark}
\begin{proof}
	We proceed by induction over the number $n$ of conjugations.
	
	For $n=1$, set $k:=k_1$, $\nu:=\lfloor P/k\rfloor$.
	Proposition~\ref{prop:single-conjugation-expansion} yields
	\begin{align*}
		A_1(t)
		=\ee^{t^kX_1}Y\ee^{-t^kX_1}
		=\sum_{j=0}^{\nu}\frac{t^{jk}}{j!}\operatorname{ad}_{X_1}^{j}(Y)+R_{\nu+1}(t),
	\end{align*}
	and since $k(\nu+1)=P-\langle P\rangle_k+k\ge P+1$, it follows that $R_{\nu+1}(t)=\mathcal{O}_{\mathcal{L}(\mathcal{D},\mathcal{H})}(t^{P+1})$.
	Define
	\begin{align*}
		B_{p,1}
		=\mathbf{1}_{k\scalebox{0.65}{$\BB{Z}$}}(p)\frac{1}{(p/k)!}\operatorname{ad}_{X_1}^{p/k}(Y),\quad
		0\le p\le P.
	\end{align*}
	Whenever $B_{p,1}\neq0$, this is a homogeneous Lie polynomial of commutator length
	\begin{align*}
		\ell(B_{p,1})
		=\frac{p}{k}\ell(X_1)+\ell(Y)
		=p+\ell(Y).
	\end{align*}
	Moreover, the highest-order commutator needed for the remainder is zero or has length
	\begin{align*}
		\ell\bigl(\operatorname{ad}_{X_1}^{\nu+1}(Y)\bigr)
		&=k(\nu+1)+\ell(Y)
		=P+\ell(Y)+k-\langle P\rangle_k\\
		&\le P+\ell(Y)+k_*
		=P^*.
	\end{align*}
	Now suppose that the assertion holds for some $n\geq1$.
	Thus,
	\begin{align}\label{eq:conjugation-induction-hypothesis}
		A_n(t)
		=\sum_{p=0}^{P}B_{p,n}t^p+R_{P+1,n}(t),
	\end{align}
	where $R_{P+1,n}(t)=\mathcal{O}_{\mathcal{L}(\mathcal{D},\mathcal{H})}(t^{P+1})$, and every $B_{p,n}$ is either zero or a homogeneous Lie polynomial of commutator length $\ell(B_{p,n})=p+\ell(Y)$.
	By definition,
	\begin{align*}
		A_{n+1}(t)
		=\ee^{t^{k_{n+1}}X_{n+1}}A_n(t)\ee^{-t^{k_{n+1}}X_{n+1}},
	\end{align*}
	and by inserting \eqref{eq:conjugation-induction-hypothesis}, we obtain
	\begin{align}\label{eq:conjugation-induction-step}
		A_{n+1}(t)
		=\sum_{p=0}^{P}t^p\ee^{t^{k_{n+1}}X_{n+1}}B_{p,n}\ee^{-t^{k_{n+1}}X_{n+1}}+\ee^{t^{k_{n+1}}X_{n+1}}R_{P+1,n}(t)\ee^{-t^{k_{n+1}}X_{n+1}}.
	\end{align}
	The last term remains of order $\mathcal{O}(t^{P+1})$ in $\mathcal{L}(\mathcal{D},\mathcal{H})$.
	Indeed, local uniform boundedness of the groups on $\mathcal{H}$ and $\mathcal{D}$ gives, for $0\le t\le T$ with fixed $T>0$,
	\begin{align*}
		&\left\|\ee^{t^{k_{n+1}}X_{n+1}}R_{P+1,n}(t)\ee^{-t^{k_{n+1}}X_{n+1}}\right\|_{\mathcal{L}(\mathcal{D},\mathcal{H})}\\
		&\qquad\leq\left\|\ee^{t^{k_{n+1}}X_{n+1}}\right\|_{\mathcal{L}(\mathcal{H})}\left\|R_{P+1,n}(t)
		\right\|_{\mathcal{L}(\mathcal{D},\mathcal{H})}\left\|\ee^{-t^{k_{n+1}}X_{n+1}}\right\|_{\mathcal{L}(\mathcal{D})}
		\le ct^{P+1},
	\end{align*}
	where $c>0$ is the product of the remainder bound from the induction hypothesis and the local group bounds on $\mathcal{H}$ and $\mathcal{D}$.
	For every $0\le p\le P$, set
	\begin{align}\label{eq:conjugation-truncation-index}
		\nu_p
		:=\left\lfloor\frac{P-p}{k_{n+1}}\right\rfloor.
	\end{align}
	Applying Proposition~\ref{prop:single-conjugation-expansion} to each component of $B_{p,n}$ and using linearity, we obtain
	\begin{align}\label{eq:coefficient-conjugation-expansion}
		\ee^{t^{k_{n+1}}X_{n+1}}B_{p,n}\ee^{-t^{k_{n+1}}X_{n+1}}
		=\sum_{j=0}^{\nu_p}\frac{t^{jk_{n+1}}}{j!}\operatorname{ad}_{X_{n+1}}^{j}(B_{p,n})+R_{\nu_p+1,p}(t),
	\end{align}
	where $R_{\nu_p+1,p}(t)=\mathcal{O}_{\mathcal{L}(\mathcal{D},\mathcal{H})}\big(t^{k_{n+1}(\nu_p+1)}\big)$.
	After multiplication by $t^p$, the corresponding remainder satisfies
	\begin{align*}
		t^pR_{\nu_p+1,p}(t)
		=\mathcal{O}_{\mathcal{L}(\mathcal{D},\mathcal{H})}\big(t^{p+k_{n+1}(\nu_p+1)}\big),
	\end{align*}
	and the definition of $\nu_p$ in \eqref{eq:conjugation-truncation-index} implies
	\begin{align*}
		p+k_{n+1}(\nu_p+1)
		=P+k_{n+1}-\langle P-p\rangle_{k_{n+1}}
		\ge P+1.
	\end{align*}
If $B_{p,n}\ne0$, the highest-order commutator required in this expansion is zero or has length
	\begin{align*}
		\ell\left(\operatorname{ad}_{X_{n+1}}^{\nu_p+1}(B_{p,n})\right)
		&=\ell(B_{p,n})+(\nu_p+1)\ell(X_{n+1})
		=P+\ell(Y)+k_{n+1}-\langle P-p\rangle_{k_{n+1}}\\[1mm]
		&\le P+\ell(Y)+k_{n+1}
		\le P+\ell(Y)+k_*
		=P^*.
	\end{align*}
	Thus, all commutators required in \eqref{eq:coefficient-conjugation-expansion} are covered by the assumptions of the theorem.
	Inserting the expansion \eqref{eq:coefficient-conjugation-expansion} into \eqref{eq:conjugation-induction-step} and collecting equal powers of $t$ gives
	\begin{align*}
		A_{n+1}(t)
		=\sum_{p=0}^{P}\sum_{j=0}^{\nu_p}\frac{t^{p+jk_{n+1}}}{j!}\operatorname{ad}_{X_{n+1}}^{j}(B_{p,n})+R_{P+1,n+1}(t)
		=\sum_{r=0}^{P}B_{r,n+1}t^r+R_{P+1,n+1}(t).
	\end{align*}
Here, the remainder $R_{P+1,n+1}(t)$ collects the conjugated remainder in \eqref{eq:conjugation-induction-step} and the new remainders from expanding the conjugations of $B_{p,n}$, each multiplied by $t^p$.
	In particular, it is $\mathcal{O}_{\mathcal{L}(\mathcal{D},\mathcal{H})}(t^{P+1})$, with a bound depending only on the local group bounds and commutators of length $\le P^*$.
	Moreover,
	\begin{align}\label{eq:conjugation-coefficient-recursion}
		B_{r,n+1}
		:=\sum_{\substack{0\le p\le r,\,j\ge 0\\[1pt]p+jk_{n+1}=r}}\frac{1}{j!}\operatorname{ad}_{X_{n+1}}^{j}(B_{p,n}).
	\end{align}
	The sum in \eqref{eq:conjugation-coefficient-recursion} is finite and every nonzero summand has commutator length
	\begin{align*}
		\ell\left(\operatorname{ad}_{X_{n+1}}^{j}(B_{p,n})\right)
		=\ell(B_{p,n})+j\ell(X_{n+1})
		=p+\ell(Y)+jk_{n+1}
		=r+\ell(Y).
	\end{align*}
	Hence $B_{r,n+1}$ is either zero or a homogeneous Lie polynomial of commutator length $r+\ell(Y)$.
\end{proof}

The finite conjugation expansions allow us to control the factored remainder $G_{2n-1}$ and derive the local-error bound under explicit regularity assumptions.

\subsection{Local error of the symmetric Zassenhaus splitting}
We now determine sufficient regularity assumptions for the local-error bound
\begin{align*}
	\|U(\tau)\psi-Z_{2n-1}(\tau)\psi\|_{\mathcal{H}}
	\le c\tau^{2n+1}\|\psi\|_{\mathcal{D}}.
\end{align*}
By the integral representation \eqref{eq:duhamel-error-representation}, it suffices to establish $N_{2n-1}(\tau)=\mathcal{O}_{\mathcal{L}(\mathcal{D},\mathcal{H})}(\tau^{2n})$.
Using the factorization $N_{2n-1}(\tau)=G_{2n-1}(\tau)Z_{2n-1}(\tau)$, we obtain this estimate by controlling the factored remainder $G_{2n-1}$ and requiring the splitting factors to preserve $\mathcal{D}$.
The bounded calculation already identifies the cancellations needed in $G_{2n-1}$: according to Theorem~\ref{thm:vanishing-defect-coefficients}, its coefficients of degrees $0,\dots,2n-1$ vanish as formal Lie polynomials.
To use these cancellations for unbounded operators, we must justify finite expansions of the conjugations in \eqref{eq:factored-defect-coefficients}, with remainders controlled in $\mathcal{L}(\mathcal{D},\mathcal{H})$.

Recall that the $\ell$-th summand in the representation of $G_{2n-1}$ in \eqref{eq:factored-defect-coefficients} carries the factor $\tau^{2\ell-2}$.
We therefore expand its conjugations up to degree
\begin{align}\label{eq:defect-truncation-degree}
	P_\ell
	:=2n-2\ell+1,\quad
	1\le\ell\le n.
\end{align}
The resulting remainder then contributes, as required,
\begin{align*}
	\tau^{2\ell-2}\mathcal{O}(\tau^{P_\ell+1})
	=\mathcal{O}(\tau^{2n}).
\end{align*}
The polynomial terms of lower degree then cancel by the simplification of algebraic commutator identities established in Theorem~\ref{thm:vanishing-defect-coefficients}.

It remains to determine which commutators these expansions require.
By \eqref{eq:simplified-factored-defect}, the conjugated polynomials are $\mathcal{R}_{2\ell-1}(1)$ for $\ell\ge2$, and $B$ for $\ell=1$ -- each of length $2\ell-1$.
Moreover, the conjugating generators are scalar multiples of $A$, $B$, and $\mathcal{R}_3(1),\dots,\mathcal{R}_{2n-1}(1)$ -- so their lengths are at most $2n-1$.
Theorem~\ref{thm:multiple-conjugation-expansion} therefore gives the uniform
sufficient commutator length bound
\begin{align}\label{eq:factored-regularity-bound}
	P_\ell+(2\ell-1)+(2n-1)
	=4n-1.
\end{align}
This leads to the following theorem, which combines control of these commutators with the flow and conjugation assumptions needed to justify the expansions and the integral error representation.
\begin{theorem}[Local error of the symmetric Zassenhaus splitting]\label{thm:unbounded-local-error}
	Let $n\ge1$, and let $\mathcal{D},\mathcal{D}_0$ be as in Assumption~\ref{ass:unbounded-framework}.
	Assume:
	\begin{enumerate}
		\item All products defining Lie polynomials in $\GG{L}_{\le 4n-1}(A,B)$ are defined on $\mathcal{D}_0$, and these polynomials extend boundedly from $\mathcal{D}$ to $\mathcal{H}$.
		\item A chosen realization of each splitting generator $\mathcal{R}_{2\ell-1}(1)$, $0\le\ell\le n$, generates a strongly continuous group on $\mathcal{H}$ that preserves $\mathcal{D}$ and is strongly continuous there.
			Moreover, a chosen realization of $A+B$ generates a strongly continuous semigroup $U(t)$ with the same properties for $t\ge 0$.
			All realizations agree with their algebraic expressions on $\mathcal{D}_0$.
		\item The integral conjugation identity in Assumption~\ref{ass:unbounded-framework} holds for every splitting generator $X$ and every $F\in\GG{L}_q(A,B)$ with $1\le q\le4n-1-\ell(X)$.
	\end{enumerate}
	Then there exist $c,\tau_0>0$ such that
	\begin{align*}
		\|U(\tau)\psi-Z_{2n-1}(\tau)\psi\|_{\mathcal{H}}
		\le c\tau^{2n+1}\|\psi\|_{\mathcal{D}},\quad
		\psi\in\mathcal{D},\quad
		0\le\tau\le\tau_0.	
	\end{align*}
\end{theorem}
\begin{proof}
	We apply Theorem~\ref{thm:multiple-conjugation-expansion} to the conjugations in the $\ell$-th summand of $G_{2n-1}$, with $P_\ell$ chosen in \eqref{eq:defect-truncation-degree}.
	By \eqref{eq:factored-regularity-bound}, the assumptions control all required commutators and yield
	\begin{align*}
		G_{2n-1}(\tau)
		=\sum_{k=0}^{2n-1}K_k\tau^k+\mathcal{O}_{\mathcal{L}(\mathcal{D},\mathcal{H})}(\tau^{2n}),
	\end{align*}
	with exactly the same Lie polynomials $K_k$ as in Lemma~\ref{lem:factored-defect-coefficients}.
	Indeed, the coefficient recursion \eqref{eq:conjugation-coefficient-recursion} coincides with the formal expansion in Proposition~\ref{prop:bounded-conjugation-expansion}, so the coefficients $B_{p,n}$ in Theorem~\ref{thm:multiple-conjugation-expansion} are the same Lie polynomials as in the bounded case.
	By Theorem~\ref{thm:vanishing-defect-coefficients}, these vanish algebraically on $\mathcal{D}_0$.
	Their bounded extensions vanish on $\mathcal{D}$ by density.
	Thus $G_{2n-1}(\tau)=\mathcal{O}_{\mathcal{L}(\mathcal{D},\mathcal{H})}(\tau^{2n})$.
	Furthermore, since the symmetric Zassenhaus splitting $Z_{2n-1}(\tau)$ is locally bounded on $\mathcal{D}$, the factorization \eqref{eq:factored-defect} gives
	\begin{align}\label{eq:unbounded-defect-order}
		N_{2n-1}(\tau)
		=\mathcal{O}_{\mathcal{L}(\mathcal{D},\mathcal{H})}(\tau^{2n}).
	\end{align}
	To justify the integral form of the local error in \eqref{eq:duhamel-error-representation}, note that each generator is closed and agrees on $\mathcal{D}_0$ with a bounded map $\mathcal{D}\to\mathcal{H}$.
	Approximation from $\mathcal{D}_0$ therefore shows that its domain contains $\mathcal{D}$.
	Invariance and strong continuity on $\mathcal{D}$ then justify the product rule in $\mathcal{H}$, giving
	\begin{align*}
		\frac{\mathrm{d}}{\mathrm{d}s}\left(U\!\left(\tfrac{\tau-s}{2}\right)Z_{2n-1}(s)U\!\left(\tfrac{\tau-s}{2}\right)\psi\right)
		=-U\!\left(\tfrac{\tau-s}{2}\right)N_{2n-1}(s)U\!\left(\tfrac{\tau-s}{2}\right)\psi.	
	\end{align*}
	Integration proves \eqref{eq:duhamel-error-representation}.
	Local semigroup bounds and \eqref{eq:unbounded-defect-order} now imply
	\begin{align*}
		\|U(\tau)\psi-Z_{2n-1}(\tau)\psi\|_{\mathcal{H}}
		\le c_0\int_0^\tau s^{2n}\,\mathrm{d}s\,\|\psi\|_{\mathcal{D}}
		=\frac{c_0}{2n+1}\tau^{2n+1}\|\psi\|_{\mathcal{D}}.	
	\end{align*}
\end{proof}

The bound $4n-1$ provides a general sufficient condition, but the regularity required in a specific application depends on which norms control the commutators that actually occur.
For differential operators, their orders may be lower than their commutator lengths suggest, allowing weaker regularity requirements, provided the flow and integral conjugation assumptions remain valid.
We conclude by examining these conditions for the Schr\"odinger equation.

\subsection{Regularity for the linear Schr\"odinger equation}
Consider $A=\ii\Delta$ and $B=-\ii V$ on $L^2(\BB{R}^d)$, where $V\in C^\infty(\BB{R}^d;\BB{R})$.
Theorem~\ref{thm:unbounded-local-error} separates the algebraic cancellations determining the local-error order from the estimates needed for a particular potential.
Once its commutator, flow, and integral conjugation assumptions are verified, no further order calculation is required.

For a general smooth potential, Schr\"odinger commutators exhibit a reduction in differential order; see \cite{Bader2014}.
More precisely, a direct induction using the Leibniz rule shows that a nonzero nested commutator of length $q\ge2$ has differential order at most $q-1$, with coefficients formed from derivatives of $V$ up to order $2q-2$.
Since the theorem requires commutators of length at most $4n-1$, their differential orders are at most $4n-2$, and their coefficients involve derivatives of $V$ up to order $2(4n-1)-2=8n-4$.
Thus, if $V$ and its derivatives up to this order are bounded, the commutator condition holds with
\begin{align*}
	\mathcal{D}
	=H^{4n-2}(\BB{R}^d),\qquad
	\mathcal{D}_0
	=C_c^\infty(\BB{R}^d).	
\end{align*}
Provided the flow and integral conjugation assumptions also hold with $H^{4n-2}(\BB{R}^d)$, the resulting local-error estimate is
\begin{align*}
	\|\ee^{\ii\tau(\Delta-V)}\psi-Z_{2n-1}(\tau)\psi\|_{L^2}
	\le c\tau^{2n+1}\|\psi\|_{H^{4n-2}}.
\end{align*}
Bounded derivatives establish the stated commutator bounds; the assumptions concerning the correction flows must be checked separately.
For unbounded potentials, weighted spaces may additionally be needed to control the coefficients.

Note that the requirement $\psi\in H^{4n-2}(\BB{R}^d)$ is sufficient under these assumptions, but need not be sharp.
It controls every commutator of length $\le 4n-1$ using a uniform differential-order estimate.
A more specific analysis can instead exploit the lower differential orders of the commutators that actually occur.
For example, Theorem~7.6 of \cite{Lasser2020} gives a fourth-order symmetric Zassenhaus error bound in terms of the semiclassically scaled $H^4$-norm, whereas the general commutator condition above gives $H^6$ when $n=2$.
The general theorem thus provides a sufficient framework, while estimates tailored to a particular equation may further reduce its regularity requirements.

The quadratic potential illustrates a different simplification: all nested commutators remain quadratic differential operators.
As a result, increasing the splitting order requires no additional regularity, and one fixed weighted $H^2$ space suffices.
\begin{example}[Quadratic Schr\"odinger potential]
	Let $A=\ii\Delta$ and $B=-\ii|x|^2$ on $L^2(\BB{R}^d)$, and set
	\begin{align*}
		\mathcal{D}
		=\{\psi\in H^2(\BB{R}^d)\,\colon\,|x|^2\psi\in L^2(\BB{R}^d)\},\qquad
		\|\psi\|_{\mathcal{D}}^2
		=\|\psi\|_{H^2}^2+\||x|^2\psi\|_{L^2}^2.		
	\end{align*}
	Then Theorem~\ref{thm:unbounded-local-error} applies for every $n\ge1$, with $\mathcal{D}_0=\mathcal{S}(\BB{R}^d)$, the Schwartz space.
	Indeed, the space
	\begin{align}\label{eq:quadratic-operator-space}
		\left\{\ii\alpha\Delta+\ii\beta|x|^2+\gamma\left(x\cdot\nabla+\tfrac d2\operatorname{Id}\right)\,\colon\,\alpha,\beta,\gamma\in\BB{R}\right\}	
	\end{align}
	is closed under commutators and contains $A$ and $B$.
	We verify the three assumptions of Theorem~\ref{thm:unbounded-local-error}.
	\begin{enumerate}
		\item Every required commutator belongs to \eqref{eq:quadratic-operator-space} and is bounded from $\mathcal{D}$ to $L^2$, using
			\begin{align*}
				\|x_r\partial_{x_q}\psi\|_{L^2}
				\le c\|\psi\|_{\mathcal{D}},\quad
				1\le r,q\le d.
			\end{align*}
			The Schwartz space $\mathcal{S}$ is dense in $\mathcal{D}$ and is preserved by differentiation and multiplication by polynomials.
			Hence all products defining these commutators are defined on $\mathcal{S}$.
		\item Each operator in \eqref{eq:quadratic-operator-space}, initially defined on $\mathcal{S}$, has a closure generating a strongly continuous unitary group on $L^2$; see e.g.\ \cite[Section~3]{Wahlberg2022}.
			By \cite[Theorem~5.2 and Corollary~5.3]{Wahlberg2022}, these groups preserve $\mathcal{S}$ and $\mathcal{D}$ and are strongly continuous on both.
		\item For a generator $X$, a nested commutator $F$, and $\psi\in\mathcal{S}$, invariance of $\mathcal{S}$ justifies differentiation on this space, giving
			\begin{align*}
				\frac{\mathrm{d}}{\mathrm{d}t}\big(\ee^{tX}F\ee^{-tX}\psi\big)
				=\ee^{tX}[X,F]\ee^{-tX}\psi.	
			\end{align*}
			Integrating from $0$ to $t$ gives the integral conjugation identity on $\mathcal{S}$.
			The boundedness of $F$ and $[X,F]$ from $\mathcal{D}$ to $L^2$, together with the local group bounds, extends this identity to $\mathcal{D}$ by density.
	\end{enumerate}
	All assumptions of Theorem~\ref{thm:unbounded-local-error} therefore hold for every $n\ge1$, and
	\begin{align*}
		\|\ee^{\ii\tau(\Delta-|x|^2)}\psi-Z_{2n-1}(\tau)\psi\|_{L^2}
		\le c_n\tau^{2n+1}\big(\|\psi\|_{H^2}^2+\||x|^2\psi\|_{L^2}^2\big)^{1/2}
	\end{align*}
	for $\psi\in\mathcal{D}$ and sufficiently small $\tau\ge0$.
	Hence arbitrarily high splitting orders are obtained on the same weighted regularity space.
\end{example}
%

\paragraph{Acknowledgments.}
The work of P.B. was funded by the Polish NCN project no.~2025/59/D/ST1/01466.
The work of K.K. was funded by the Polish NCN project no.~2019/34/E/ST1/00390.
We thank Caroline Lasser for helpful discussions on domain questions for unbounded operators and for drawing our attention to weighted Sobolev spaces.


\bibliographystyle{agsm}
\bibliography{Zassenhaus}

\appendix
\renewcommand{\theequation}{\Alph{section}.\arabic{equation}}
\setcounter{equation}{0}
\section{Detailed derivation of the Zassenhaus exponents for $Z_9$}\label{app:zassenhaus-exponents-z9}
In what follows, we will present a detailed derivation of the symmetric Zassenhaus splitting $Z_9$.
Note that this corresponds to the choice $n=5$.
Let us briefly recall the definition of the arguments in this specific case:
\begin{align*}
	\mathcal{B}_{-1}
	&:=(-\tau A,\tau A+\tau B),\\
	\mathcal{B}_1
	&:=(-\tau B,\mathrm{sBCH}_9(\mathcal{B}_{-1})),\\
	\mathcal{B}_{2k-1}
	&:=(-2\mathcal{R}_{2k-1},\mathrm{sBCH}_9(\mathcal{B}_{2k-3})),
	\quad 2\le k\le 5,
\end{align*}
and the following index sets for $2\le\ell\le 5$:
\begin{align*}
	\mathcal{I}_{9,2\ell-1}
	&:=\{(j,k)\in\BB{Z}^2\,\colon\,2\le j\le 5,\,0\le k\le\ell-1\},\\[2mm]
	\mathcal{J}^{2\ell-1}
	&:=\{(j,k)\in\BB{Z}^2\,\colon\,j\ge2,\,k\ge0,\,(2j-1)|2k-1|\le2\ell-1\}.
\end{align*}
\textbf{Step 1 ($\ell=1$):}
\begin{align*}
	\ee^{-\frac12\tau A}\ee^{\tau A+\tau B}\ee^{-\frac12\tau A}
	&=\ee^{\mathrm{sBCH}_9({\mathcal{B}}_{-1})}+E_{9}^{\operatorname{sBCH}}({\mathcal{B}}_{-1})\\
	&=\ee^{\tau B+\sum_{j=2}^5C_{2j-1}({\mathcal{B}}_{-1})}+E_{9}^{\operatorname{sBCH}}({\mathcal{B}}_{-1}).
\end{align*}
Leading term of $\mathrm{sBCH}_9({\mathcal{B}}_{-1})$: $\tau B$.
Error: $E_{9}^{\operatorname{sBCH}}({\mathcal{B}}_{-1})=\mathcal{O}(\tau^{1\times 11})=\mathcal{O}(\tau^{11})$.

\noindent
\textbf{Step 2 ($\ell=2$):}
The index sets are given by
\begin{align*}
	\mathcal{I}_{9,3}
	&=\{(2,0),(2,1),(3,0),(3,1),(4,0),(4,1),(5,0),(5,1)\},\\
	\mathcal{J}^{3}
	&=\{(2,0),(2,1)\}.
\end{align*}
\begin{align*}
	\ee^{-\frac12\tau B}\ee^{\mathrm{sBCH}_9({\mathcal{B}}_{-1})}\ee^{-\frac12\tau B}
	&=\ee^{\mathrm{sBCH}_9({\mathcal{B}}_1)}+E_{9}^{\operatorname{sBCH}}({\mathcal{B}}_1)\\
	&=\ee^{-\tau B+\mathrm{sBCH}_9({\mathcal{B}}_{-1})+\sum_{j=2}^5C_{2j-1}({\mathcal{B}}_1)}+E_{9}^{\operatorname{sBCH}}({\mathcal{B}}_1)\\
	&=\ee^{2\mathcal{R}_3(\tau)+{\mathcal{O}(\tau^{5})}}+E_{9}^{\operatorname{sBCH}}({\mathcal{B}}_1).
\end{align*}
Leading term of $\mathrm{sBCH}_9({\mathcal{B}}_{1})$:
\begin{align*}
	2\mathcal{R}_3(\tau)
	=\tau^3[\tau^3]\sum_{(j,k)\in\mathcal{J}^3}C_{2j-1}(\mathcal{B}_{2k-1})
	=\tau^3[\tau^3]\left(C_3(\mathcal{B}_{-1})+C_3(\mathcal{B}_1)\right).
\end{align*}
Error: $E_{9}^{\operatorname{sBCH}}({\mathcal{B}}_1)=\mathcal{O}(\tau^{1\times 11})=\mathcal{O}(\tau^{11})$.

\noindent
\textbf{Step 3 ($\ell=3$):}
The index sets are given by
\begin{align*}
	\mathcal{I}_{9,5}
	&=\{(2,0),(2,1),(2,2),(3,0),\dots,(5,2)\},\\
	\mathcal{J}^{5}
	&={\{(2,0),(2,1),(3,0),(3,1)\}}.
\end{align*}
\begin{align*}
	\ee^{-\mathcal{R}_3(\tau)}\ee^{\mathrm{sBCH}_9({\mathcal{B}}_1)}\ee^{-\mathcal{R}_3(\tau)}
	&=\ee^{\mathrm{sBCH}_9({\mathcal{B}}_3)}+E_{9}^{\operatorname{sBCH}}({\mathcal{B}}_3)\\
	&=\ee^{-2\mathcal{R}_3(\tau)+\mathrm{sBCH}_9({\mathcal{B}}_1)+\sum_{j=2}^5 C_{2j-1}({\mathcal{B}}_3)}+E_{9}^{\operatorname{sBCH}}({\mathcal{B}}_3)\\
	&=\ee^{2\mathcal{R}_5(\tau)+{\mathcal{O}(\tau^{7})}}+E_{9}^{\operatorname{sBCH}}({\mathcal{B}}_3).
\end{align*}
Leading term of $\mathrm{sBCH}_9({\mathcal{B}}_{3})$:
\begin{align*}
	2\mathcal{R}_5(\tau)
	=\tau^5[\tau^5]\sum_{(j,k)\in\mathcal{J}^5}C_{2j-1}(\mathcal{B}_{2k-1})
	=\tau^5[\tau^5]\left(\sum_{j=2}^{3}\bigl(C_{2j-1}(\mathcal{B}_{-1})+C_{2j-1}(\mathcal{B}_1)\bigr)\right).
\end{align*}
Error: $E_{9}^{\operatorname{sBCH}}({\mathcal{B}}_3)=\mathcal{O}(\tau^{3\times 11})=\mathcal{O}(\tau^{33})$.

\noindent
\textbf{Step 4 ($\ell=4$):}
The index sets are given by
\begin{align*}
	\mathcal{I}_{9,7}
	&=\{(2,0),(2,1),(2,2),(2,3),(3,0),\dots,(5,3)\},\\
	\mathcal{J}^{7}
	&={\{(2,0),(2,1),(3,0),(3,1),(4,0),(4,1)\}}.
\end{align*}
\begin{align*}
	\ee^{-\mathcal{R}_5(\tau)}\ee^{\mathrm{sBCH}_9({\mathcal{B}}_3)}\ee^{-\mathcal{R}_5(\tau)}
	&=\ee^{\mathrm{sBCH}_9({\mathcal{B}}_5)}+E_{9}^{\operatorname{sBCH}}({\mathcal{B}}_5)\\
	&=\ee^{-2\mathcal{R}_5(\tau)+\mathrm{sBCH}_9({\mathcal{B}}_3)+\sum_{j=2}^5 C_{2j-1}({\mathcal{B}}_5)}+E_{9}^{\operatorname{sBCH}}({\mathcal{B}}_5)\\
	&=\ee^{2\mathcal{R}_7(\tau)+{\mathcal{O}(\tau^{9})}}+E_{9}^{\operatorname{sBCH}}({\mathcal{B}}_5).
\end{align*}
Leading term of $\mathrm{sBCH}_9({\mathcal{B}}_{5})$:
\begin{align*}
	2\mathcal{R}_7(\tau)
	=\tau^7[\tau^7]\sum_{(j,k)\in\mathcal{J}^7}C_{2j-1}(\mathcal{B}_{2k-1})
	=\tau^7[\tau^7]\left(\sum_{j=2}^{4}\bigl(C_{2j-1}(\mathcal{B}_{-1})+C_{2j-1}(\mathcal{B}_1)\bigr)\right).
\end{align*}
Error: $E_{9}^{\operatorname{sBCH}}({\mathcal{B}}_5)=\mathcal{O}(\tau^{5\times 11})=\mathcal{O}(\tau^{55})$.

\noindent
\textbf{Step 5 ($\ell=5$):}
The index sets are given by
\begin{align*}
	\mathcal{I}_{9,9}
	&=\{(2,0),(2,1),(2,2),(2,3),(2,4),(3,0),\dots,(5,4)\},\\
	\mathcal{J}^{9}
	&={\bigl(\{2,3,4,5\}\times\{0,1\}\bigr)\cup\{(2,2)\}}.
\end{align*}
\begin{align*}
	\ee^{-\mathcal{R}_7(\tau)}\ee^{\mathrm{sBCH}_9({\mathcal{B}}_5)}\ee^{-\mathcal{R}_7(\tau)}
	&=\ee^{\mathrm{sBCH}_9({\mathcal{B}}_7)}+E_{9}^{\operatorname{sBCH}}({\mathcal{B}}_7)\\
	&=\ee^{-2\mathcal{R}_7(\tau)+\mathrm{sBCH}_9({\mathcal{B}}_5)+\sum_{j=2}^5 C_{2j-1}({\mathcal{B}}_7)}+E_{9}^{\operatorname{sBCH}}({\mathcal{B}}_7)\\
	&=\ee^{2\mathcal{R}_9(\tau)+{\mathcal{O}(\tau^{11})}}+E_{9}^{\operatorname{sBCH}}({\mathcal{B}}_7).
\end{align*}
Leading term of $\mathrm{sBCH}_9({\mathcal{B}}_{7})$:
\begin{align*}
	2\mathcal{R}_9(\tau)
	=\tau^9[\tau^9]\sum_{(j,k)\in\mathcal{J}^9}C_{2j-1}(\mathcal{B}_{2k-1})
	=\tau^9[\tau^9]\left(\sum_{j=2}^{5}\bigl(C_{2j-1}(\mathcal{B}_{-1})+C_{2j-1}(\mathcal{B}_1)\bigr)+C_3(\mathcal{B}_3)\right).
\end{align*}
Error: $E_{9}^{\operatorname{sBCH}}({\mathcal{B}}_7)=\mathcal{O}(\tau^{7\times 11})=\mathcal{O}(\tau^{77})$.

\setcounter{equation}{0}
\section{Derivatives of the conjugation map $t\mapsto\ee^{t^kX}Y\ee^{-t^kX}$}\label{app:conjugation-derivatives}
For the computation of the derivatives $H^{(s)}$, $1\le s\le\nu+1$, we let
\begin{align*}
	F(t)
	=\ee^{tX}Y\ee^{-tX}
	\quad\text{and}\qquad
	g(t)
	=t^k,
\end{align*}
with derivatives given by
\begin{align*}
	F^{(j)}(t)
	=\ee^{tX}\operatorname{ad}_X^j(Y)\ee^{-tX},\quad
	g^{(j)}(t)
	=\begin{cases}[k]_jt^{k-j},&1\le j\le k,\\0,&j>k,\end{cases},\quad
	j\ge 1,
\end{align*}
where $[k]_j=k(k-1)\cdots(k-j+1)$ denotes the descending factorial.
Recall that $[k]_j=0$ if $j>k$ by definition.
Then, we have $H=F\circ g$ and applying Fa\`{a} di Bruno's formula, we obtain that the $s$-th derivative can be computed as
\begin{align}\label{eq:faa-di-bruno}
	H^{(s)}(t)
	=\sum_{m\in\mathcal{M}_s}d_{m,s}F^{(|m|)}(g(t))\prod_{j=1}^s\left(g^{(j)}(t)\right)^{m_j},\quad
	0^0
	:=1,
\end{align}
where $\mathcal{M}_s=\{m\in\BB{Z}_{\ge 0}^s\,\colon\,\sum_{j=1}^sjm_j=s\}$, and the constants $d_{m,s}$ are given by
\begin{align*}
	d_{m,s}
	:=\frac{s!}{\prod_{j=1}^sm_j!(j!)^{m_j}}.
\end{align*}
Substituting the above expressions for the derivatives of $g$ gives, for $s>k$ and $t\ne0$,
\begin{align*}
	\prod_{j=1}^s\left(g^{(j)}(t)\right)^{m_j}
	=\prod_{j=1}^{k}[k]_j^{m_j}t^{m_j(k-j)}\prod_{j=k+1}^{s}0^{m_j}
	=t^{k|m|-s}\prod_{j=1}^{k}[k]_j^{m_j}\prod_{j=k+1}^{s}0^{m_j}.
\end{align*}
Since the product $\prod_{j=k+1}^{s}0^{m_j}$ is nonzero if and only if $m_j=0$ for all $j>k$, it suffices to replace the sum over $\mathcal{M}_s$ in \eqref{eq:faa-di-bruno} by a sum over the smaller index set $\mathcal{M}_{k,s}$ defined in \eqref{eq:derivative-multiindices}.
Note that in particular $\mathcal{M}_{k,s}=\mathcal{M}_s$ if $s\le k$.
Substituting the above expressions for the derivatives of $F$ gives
\begin{align*}
	H^{(s)}(t)
	&=\sum_{m\in\mathcal{M}_{k,s}}d_{m,s}\ee^{t^kX}\operatorname{ad}_X^{|m|}(Y)\ee^{-t^kX}t^{k|m|-s}\prod_{j=1}^{s}[k]_j^{m_j}\\
	&=\sum_{m\in\mathcal{M}_{k,s}}c_{m,s}\ee^{t^kX}\operatorname{ad}_X^{|m|}(Y)\ee^{-t^kX}t^{k|m|-s}
\end{align*}
The resulting formula also holds at $t=0$ by continuity, with the constants
\begin{align}\label{eq:derivative-combinatorial-coefficients}
	c_{m,s}
	:=d_{m,s}\prod_{j=1}^{s}[k]_j^{m_j},\quad
	m\in\mathcal{M}_{k,s}.
\end{align}
Using that $\ee^{t^kX}\big|_{t=0}=\operatorname{Id}$ and that $0^{k|m|-s}\ne 0$ only if $s$ is a multiple of $k$, we see that $H^{(s)}(0)=0$ if $s\bmod k\ne0$.\footnote{For a given $s=jk$, choose $m=j\mathbf{e}_k\in\mathcal{M}_{k,s}$. Then $|m|=j$.}
If $s\bmod k=0$, set $j=s/k$.
Then
\begin{align*}
	H^{(s)}(0)
	=\frac{s!}{j!(k!)^j}(k!)^j\operatorname{ad}_X^j(Y)
	=\frac{s!}{j!}\operatorname{ad}_X^j(Y).
\end{align*}
%

\setcounter{equation}{0}
\section{Proof of the combinatorial identity}\label{app:combinatorial-identity}
The proof of the combinatorial identity \eqref{eq:combinatorial-identity} proceeds in two stages.
First, we will show that
\begin{align}\label{eq:combinatorial-coefficient-formula}
	s_{j,\nu+1}
	=\frac{(\nu+1)!}{j!}[x^{\nu+1}]\left[(1+x)^k-1\right]^j,
\end{align}
where $[x^{\nu+1}]P(x)$ denotes the $(\nu+1)$-th coefficient of the polynomial $P$.
Second, we use the expansion of
\begin{align*}
	\left[(1+x)^k-1\right]^j
\end{align*}
to read off the corresponding coefficient directly, thereby establishing the combinatorial identity.
We start by proving \eqref{eq:combinatorial-coefficient-formula}.
We observe that for fixed $k\ge 1$ and $x\in\BB{R}$ it holds that
\begin{align*}
	(1+x)^k-1
	=\sum_{r=1}^k\binom{k}{r}x^r
\end{align*}
We regard each $m\in\mathcal{M}_{k,\nu+1}$ as having $k$ components by adding zero components if $k>\nu+1$ and omitting the final components, which are zero by definition, if $k<\nu+1$.
Hence, upon introducing the constant
\begin{align*}
	a_{k,m}:=\prod_{r=1}^k\frac{1}{m_r!}\binom{k}{r}^{m_r},\quad
	m\in\BB{Z}_{\ge0}^{k},\quad k\ge1,
\end{align*}
we can write
\begin{align*}
	s_{j,\nu+1}
	=(\nu+1)!\sum_{\begin{subarray}{c}m\in\mathcal{M}_{k,\nu+1}\\|m|=j\end{subarray}}a_{k,m}.
\end{align*}
The multinomial expansion yields
\begin{align*}
	\frac{(\nu+1)!}{j!}[x^{\nu+1}]\left[(1+x)^k-1\right]^j
	&=\frac{(\nu+1)!}{j!}[x^{\nu+1}]\left(\sum_{r=1}^k\binom{k}{r}x^r\right)^j\\
	&=\frac{(\nu+1)!}{j!}[x^{\nu+1}]\left(\sum_{m_1+\dots+ m_k=j}j!a_{k,m}\,x^{\sum_{r=1}^krm_r}\right).
\end{align*}
The $(\nu+1)$-th coefficient is given by all the summands for which we have $\sum_{r=1}^krm_r=\nu+1$, which by definition is exactly the sum over the set
\begin{align*}
	\mathcal{M}_{k,\nu+1}\cap\left\{m\in\BB{Z}_{\ge0}^k\,\colon\,|m|=j\right\},
\end{align*}
and we conclude that
\begin{align*}
	\frac{(\nu+1)!}{j!}[x^{\nu+1}]\left[(1+x)^k-1\right]^j
	&=(\nu+1)!\sum_{\begin{subarray}{c} m\in\mathcal{M}_{k,\nu+1}\\|m|=j\end{subarray}}a_{k,m}
	=s_{j,\nu+1}.
\end{align*}
For the second part of the proof, we observe that
\begin{align*}
	\sum_{j=1}^{r}\frac{s_{j,\nu+1}}{(r-j)!}
	&=\sum_{j=1}^{r}\frac{(\nu+1)!}{j!(r-j)!}[x^{\nu+1}]\left[(1+x)^k-1\right]^j\\
	&=\frac{(\nu+1)!}{r!}[x^{\nu+1}]\sum_{j=1}^{r}\binom{r}{j}\left[(1+x)^k-1\right]^j
\end{align*}
Using that we can add the term for $j=0$ since $[x^{\nu+1}]1=0$, we further conclude that
\begin{align*}
	\sum_{j=1}^{r}\frac{s_{j,\nu+1}}{(r-j)!}
	&=\frac{(\nu+1)!}{r!}[x^{\nu+1}]\sum_{j=0}^{r}\binom{r}{j}\left[(1+x)^k-1\right]^j1^{r-j}\\
	&=\frac{(\nu+1)!}{r!}[x^{\nu+1}](1+x)^{kr}
	=\frac{(\nu+1)!}{r!}\binom{kr}{\nu+1}
	=\frac{[kr]_{\nu+1}}{r!},
\end{align*}
which is exactly the combinatorial identity in \eqref{eq:combinatorial-identity}.
\end{document}